\documentclass[preprint,12pt]{elsarticle}

\usepackage[utf8]{inputenc}
\DeclareUnicodeCharacter{211D}{$\mathbb{R}$}
\usepackage{amsmath}
\usepackage{amssymb}
\usepackage{xcolor}
\usepackage{booktabs}
\usepackage{subcaption}
\usepackage{microtype}
\usepackage{lineno}
\usepackage{hyperref}
\usepackage{ulem}
\usepackage{orcidlink}

\newcommand{\resultgraphic}[2][]{%
    \IfFileExists{#2}{\includegraphics[#1]{#2}}{%
        \fbox{\parbox[c][3cm][c]{0.9\linewidth}{\centering\footnotesize
            \textcolor{red}{TODO: missing figure}\\ \nolinkurl{#2}}}}}

\newtheorem{lemma}{Lemma}
\newtheorem{corollary}{Corollary}

\journal{Journal of Computational Physics}

\begin{document}

\begin{frontmatter}

\title{Factorized approximate inverse preconditioners in the tensor train format}

\affiliation[label1]{organization={Barcelona Supercomputing Center},
            city={Barcelona},
            country={Spain}
}

\affiliation[label2]{organization={Universitat Politècnica de Catalunya},
            city={Barcelona},
            country={Spain}
}

\affiliation[label3]{organization={Qilimanjaro Quantum Tech},
            city={Barcelona},
            country={Spain}
}

\author[label1,label2]{Pau Torrente Badia\orcidlink{0009-0005-2942-3263}\corref{cor1}}
\cortext[cor1]{\href{mailto:pau.torrente@bsc.es}{pau.torrente@bsc.es}}

\author[label1]{Daniel Pérez García\orcidlink{0009-0009-0283-7019}}
\author[label1,label3]{Artur Garcia-Saez\orcidlink{0000-0003-3561-0223}}
\author[label1]{Stefano Carignano\orcidlink{0000-0001-6549-5364}}

\begin{abstract}
We present a preconditioning pipeline that accelerates Tensor Train (TT) linear solvers for symmetric positive-definite problems. 
Extending the Factorized Sparse Approximate Inverse framework to tensor networks, we construct a preconditioner 
as the product of a TT operator and its transpose. This formulation is positive semidefinite for any bond dimension and 
truncation, and definite whenever the factor is non-singular, satisfying the mathematical requirements of the 
Conjugate Gradient method. We detail the construction of this preconditioner and its integration into Density Matrix Renormalization Group 
(DMRG)-like linear solvers. Performance evaluations show consistent iteration reductions for anisotropic diffusion and Poisson problems across various configurations, leading
to net speedups of up to 25\% in ill-conditioned settings of time stepping schemes that can reuse the 
approximate inverses. These benchmarks were performed with the open-source Julia 
library DifferenTT.jl, which we also introduce in this work to provide out-of-the-box access to these preconditioned pipelines.
\end{abstract}


\end{frontmatter}


\section{Introduction}

Tensor networks provide compact representations of high-dimensional objects by decomposing them into low-rank factors,
making it possible to control computational costs while retaining the essential structure of the underlying problem \cite{Grasedyck2013,orus2014}. Among
these formats, the Tensor Train (TT) representation \cite{oseledets2011} has become a particularly useful tool for large-scale numerical
modeling \cite{garcia-ripoll2021}. Introduced by the quantum physics community under the name of Matrix Product States (MPS) \cite{schollwock2011}, these 
array chains have become one of the main tools for probing the properties of quantum systems on lattices, especially in 
one-dimensional settings \cite{klumper_1993,vidal2003}. Even though more complex networks have been introduced for systems with intricate correlation structures,
the ease of handling these chains and their computational efficiency make them exceptionally reliable \cite{orus2014}.

The mathematical community later rediscovered these structures as the TT format \cite{oseledets2011}, demonstrating their capabilities
as a generic numerical tool for high-dimensional problems such as function interpolation \cite{savostyanov2011a}, integration \cite{dolgov2020}, 
and the solution of Partial Differential Equations (PDEs) \cite{bachmayr2023}. 

A large fraction of the computational cost in these TT-based PDE solvers lies in solving the resulting large-scale linear systems ($A\mathbf{x} = \mathbf{b}$). 
Within the context of TTs, this is typically 
approached either by performing a global optimization in the network manifold \cite{dolgov2013}, 
or through alternating local updates in sweeps \cite{oseledets2012,holtz2012, dolgov2014}, analogous to the Density Matrix 
Renormalization Group (DMRG) algorithm \cite{schollwock2011}. In both approaches, iterative methods are frequently 
invoked to solve the underlying linear system \cite{saad2003,Hackbusch2016}. In particular, for symmetric positive-definite 
(SPD) linear systems, the Conjugate Gradient (CG) method \cite{Hestenes1952,barrett1994} is usually the solver of choice, which 
proceeds by minimizing the quadratic cost function
\begin{equation}
    f(\mathbf{x}) = \frac{1}{2}\mathbf{x}^T A \mathbf{x} - \mathbf{x}^T\mathbf{b}
\end{equation}
by generating a sequence of mutually $A$-conjugate search directions. 

The efficiency of the CG method, and related iterative methods, depends on the spectral properties of the operator $A$ \cite{Hackbusch2016}. In particular, the $A$-norm error of the CG iterations,
starting from an initial guess $\mathbf{x}_0$ is bounded by
\begin{equation}
    \label{eq:CG convergence}
    \|\mathbf{e}_k\|_A \le 2 \left( \frac{\sqrt{\kappa(A)} - 1}{\sqrt{\kappa(A)} + 1} \right)^k \|\mathbf{e}_0\|_A,
\end{equation}
where $\mathbf{e}_k = \mathbf{x}_k - \mathbf{x}^*$ is the error at iteration $k$ with respect to the exact solution $\mathbf{x}^*$, $\|\mathbf{v}\|_A = \sqrt{\mathbf{v}^TA\mathbf{v}}$ the $A$-norm 
and $\kappa(A)$ the condition number of the operator. Since the operator is SPD, this quantity
is simply the ratio between its largest and smallest eigenvalues. A completely clustered spectrum, such as that of a scaled identity,
lets CG converge in a single iteration. Conversely, the more spread out the spectrum is,
the more sensitive the iteration becomes to numerical errors, and convergence slows down significantly. 

In practice, unfortunately, this condition number often scales unfavorably with the system size or grid resolution. 
This deterioration is typically mitigated by means of preconditioning techniques \cite{saad2003}, which transform the linear system in order to cluster 
the eigenvalues and accelerate convergence. In this work, we focus on preconditioning techniques based on approximate inverses $M^{-1} \approx A^{-1}$, which 
modify the linear system as
\begin{equation}
    M^{-1}A\mathbf{x} = M^{-1}\mathbf{b},
\end{equation}
such that the effective condition number of the joint operator $M^{-1}A$ is reduced. A good preconditioner, therefore, must
be a fairly accurate approximation of the exact inverse, while simultaneously remaining computationally cheap to compute and apply. 
Crucially, for the preconditioned CG (PCG) method, the approximate inverse $M^{-1}$ must be SPD, just like $A$.

While SPD-preserving preconditioning is mature in classical sparse linear algebra \cite{janna2015}, not many developments have been made from the point of view
of tensor networks. Efficient and exact representations of relevant operators in the TT format only exist in some very
particular settings \cite{Kazeev2012,khoromskij2009b}, and approximate inversion techniques \cite{oseledets2010,oseledets2012} do not guarantee positive definiteness, which can make PCG break down. 

In this work we address this issue by developing an approximate inversion pipeline for general SPD operators that satisfies PCG's
requirements by adapting the Factorized Sparse Approximate Inverse (FSAI) pipeline \cite{kolotilina1993} to
the TT format. Rather than approximating the inverse of the full operator directly, this approach constructs an approximate 
inverse of its Cholesky factor. This factorized construction preserves the SPD property required by conjugate 
gradient methods by design, given that the factors are non-singular. We begin by introducing TTs and TT-based linear solvers in Sec.~\ref{sec:tensor_trains}.
We then show how to adapt the full FSAI pipeline to this format in Sec.~\ref{sec:ttfai}. By means of a projection step, 
we demonstrate how the construction can be obtained with an efficient tensor network representation.
This approximate inverse can then be applied to TT-based solvers. We focus on DMRG-based approaches 
for their robustness, and discuss in Sec.~\ref{sec:Preconditioning ALS} 
how to incorporate the factorized preconditioner into this pipeline in a way that keeps their properties without unnecessary costs.
Finally, in order to showcase the performance of this TT-based factorized approximate inverse (TT-FAI)
we evaluate it in a high-dimensional anisotropic diffusion problem and on the pressure Poisson
equation of an incompressible flow solver. These methods are implemented in DifferenTT.jl \cite{DifferenTT}, an open-source Julia library built on the ITensors.jl \cite{ITensors} ecosystem, which we introduce in this work as a workbench for solving PDEs in the TT format. These numerical experiments not only show that
the preconditioner is able to reduce the required iterations, as intended, but also that its construction time can be amortized in time stepping schemes that reuse all the operators. 

\section{Tensor trains}
\label{sec:tensor_trains}

We begin by briefly reviewing the TT format in order to establish the mathematical basis upon which the methods in this 
work are built. As introduced earlier, the primary motivation for tensor networks is to circumvent the 
\textit{curse of dimensionality} when dealing with high-dimensional problems. 

A generic order-$N$ tensor (containing the values of a function evaluated on an $N$-dimensional grid, for example) whose local degree of freedom takes $d$ values, requires $\mathcal{O}(d^N)$ entries to be stored. 
The TT format can drastically reduce this memory footprint by expressing the components of such a high-dimensional tensor $\mathcal{V}_{i_1, i_2, \dots, i_N}$ as a sequence of products
\begin{equation}
    \mathcal{V}_{i_1, i_2, \dots, i_N} = \sum_{\alpha_0, \dots, \alpha_N} V^{(1)}_{\alpha_0, i_1, \alpha_1} V^{(2)}_{\alpha_1, i_2, \alpha_2} \dots V^{(N)}_{\alpha_{N-1}, i_N, \alpha_N},
    \label{eq:mps}
\end{equation}
where $i_k \in \{1, \dots, d\}$ are the \textit{physical} indices and $\alpha_k \in \{1, \dots, \chi_k\}$ are the 
\textit{bond} indices, of dimension $\chi_k$, with the boundary conditions $\alpha_0 = \alpha_N = 1$. The maximum bond dimension, 
$\chi = \max_k \chi_k$, bounds the amount of correlations between each bipartition of the system. If $\mathcal{V}_{i_1, i_2, \dots, i_N}$
can be represented by a bounded $\chi$, then the storage cost drops to $\mathcal{O}(N d \chi^2)$, avoiding exponential growth with $N$.

Similarly, a high-dimensional linear operator $\mathcal{W}$, which in general has $\mathcal{O}(d^{2N})$ 
entries, can be factorized as a TT operator (or Matrix Product Operator, MPO). Its components are written as
\begin{equation}
    \mathcal{W}_{i_1, \dots, i_N; j_1, \dots, j_N} = \sum_{\beta_0, \dots, \beta_N} W^{(1)}_{\beta_0, i_1, j_1, \beta_1} W^{(2)}_{\beta_1, i_2, j_2, \beta_2} \dots W^{(N)}_{\beta_{N-1}, i_N, j_N, \beta_N},
    \label{eq:mpo}
\end{equation}
where $i_k$ and $j_k$ represent the physical indices. Viewing each core $W^{(k)}$ as a
$\chi_{k-1}\times\chi_k$ block matrix with $d_k\times d_k$ blocks, Eq.~\eqref{eq:mpo} reads
$\mathcal{W} = W^{(1)}\bowtie\dots\bowtie W^{(N)}$, where $\bowtie$ is the strong Kronecker product
\cite{Kazeev2012}, $(A \bowtie B)_{\alpha\gamma} = \sum_\beta A_{\alpha\beta}\otimes B_{\beta\gamma}$.

These decompositions can also be applied to one-dimensional vectors and, by extension, to the operators that act on them.
By means of a binary encoding of its index, through the expansion
\begin{equation}
    i \leftrightarrow (i_1,\dots,i_N), \;\; i = \sum_{p = 1}^N i_p 2^{N-p},  i_p \in \{0, 1\}
\end{equation}
we can reshape a $2^N$-dimensional vector into a tensor with $N$ indices of size $2$ and apply the TT
decomposition. This construction, commonly known as Quantized Tensor Train (QTT) \cite{khoromskij2011,oseledets2010}, can provide an 
exponential compression, especially for problems with uncoupled scales \cite{Gourianov2022}.

\begin{figure}[tpb]
    \centering
    \begin{subfigure}[b]{0.45\textwidth}
        \centering
        \includegraphics[width=\textwidth]{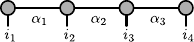} 
        \caption{Matrix Product State (MPS)}
        \label{fig:mps}
    \end{subfigure}\hfill
    \begin{subfigure}[b]{0.45\textwidth}
        \centering
        \includegraphics[width=\textwidth]{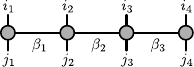}
        \caption{Matrix Product Operator (MPO)}
        \label{fig:mpo_diagram}
    \end{subfigure}
    \caption{Diagrammatic representation of (a) a TT / MPS (Eq.~\eqref{eq:mps}) and (b) a TT operator / MPO (Eq.~\eqref{eq:mpo}). The horizontal lines correspond to the contracted bond indices ($\alpha_k, \beta_k$) bounded by the bond dimension $\chi$, while vertical lines denote the free physical indices ($i_k, j_k$).}
    \label{fig:mps_mpo}
\end{figure}

These decompositions are most intuitively visualized using diagrammatic tensor notation
(see Fig.~\ref{fig:mps_mpo}). In this graphical language, each tensor is represented 
by a shape (e.g., a circle or a square), and its indices are depicted as legs. Connected legs between two tensors imply 
a summation (contraction) over that shared index. The TT representation of a vector (Eq.~\eqref{eq:mps}), thus, resembles 
a chain of tensors, each with one free physical leg, while a TT operator (Eq.~\eqref{eq:mpo}) has two free physical legs per 
site.

\subsection{Linear solvers in the TT format}
\label{subsec:tt_linear_solvers}
Having defined the TT format, we now briefly introduce how these factorized structures can be leveraged to solve 
large-scale linear systems. In the context of TTs, linear solvers can be broadly categorized into two main classes. 

The first class comprises global solvers, conceptually very similar to classical iterative solvers. They optimize the tensor 
train as a whole, leveraging implementations of general linear algebra operations in the TT format \cite{oseledets2011}, 
to perform the iterations of classical Krylov subspace methods directly within the restricted TT manifold \cite{kressner2010,ballani2013,dolgov2013}. To keep 
memory and computational costs under control, the bond dimensions of the network are periodically truncated during 
the iterative process. It is worth noting that these truncation steps do introduce errors that can break the 
convergence guarantees of some algorithms. In particular, for the CG method, severe truncation can destroy the orthogonality
between search directions, making the method fragile in practice. For this reason,
TT solvers of this kind are not considered in this work.

The second prominent approach for solving linear systems is closely related to the DMRG framework, originally developed for computing the properties of quantum systems \cite{white1992,schollwock2011}. 
Instead of optimizing the TT globally, DMRG-based solvers proceed in an alternating fashion along the network 
\cite{oseledets2012,holtz2012,dolgov2014}. This is achieved by utilizing the current solution guess to project the global 
high-dimensional linear system into a tractable local problem, focused on a small subpart of the network (typically one 
or two adjacent sites). Within this constrained local subspace, one invokes a classical iterative solver. Because the 
local optimization is exact within the current subspace, the convergence properties of this inner solver are unaffected 
by truncation errors during its iterations. Once the local solve converges, the updated tensor is reinserted into the 
tensor network. During this step, one can perform a truncated Singular Value Decomposition (SVD) to bound the bond dimension.

\begin{figure}[tpb]
    \centering
    \begin{subfigure}[c]{0.4\textwidth}
        \centering
        \includegraphics[width=\textwidth]{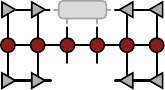}
        \caption{Projected operator $A_{\text{loc}} = U^T AU$}
        \label{fig:projectedA}
    \end{subfigure}\hspace{0.1\textwidth}
    \begin{subfigure}[c]{0.4\textwidth}
        \centering
        \includegraphics[width=\textwidth]{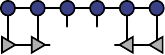}
        \caption{Projected right-hand side $b_{\text{loc}} = U^T b$}
        \label{fig:projectedb}
    \end{subfigure}
    \caption{Projected operator and right-hand side involved in a two-site MALS pipeline where the 
 optimization is taking place on the third and fourth sites of a TT.
     The dark red tensors correspond to the operator $A$
    and the blue tensors to $b$. In both of them, the gray triangles (pointing towards the orthogonality center) that form the projector 
    correspond to the isometric tensors derived from the current solution TT. The blurred rectangular gray tensor depicts the local tensor $x_{\text{loc}}$ that is being
    updated by the local solver to reach a local solution.}
    \label{fig:intro ALS}
\end{figure}

This framework offers significant flexibility in how the projection and subsequent recovery of the tensor network 
structure are handled. For instance, projecting the problem onto multiple contiguous sites simultaneously can help 
avoid local minima and allows for dynamic adjustment of the bond dimension \cite{holtz2012}, while basis expansion 
techniques can be employed to guide the optimization \cite{dolgov2014}. For illustrative purposes, 
Fig.~\ref{fig:intro ALS} depicts the projection process for the two-site variant of the Alternating Linear Scheme (ALS) \cite{holtz2012}, known as Modified ALS (MALS).
In the figure, the triangular sites that form 
the projection are isometries derived from the canonical form of the current solution guess \cite{schollwock2011}. Once the local
problem is solved, the TT format is recovered by means of an SVD and the
optimization proceeds onto the next pair of sites, continuing with the sweeping procedure. Crucially, as we will 
detail in Sec.~\ref{sec:Preconditioning ALS}, this projection mechanism applies not only to the main operator $A$ of 
the linear system, but can also be used to systematically apply approximate inverses to precondition the local systems 
\cite{kressner2011}.

\section{Factorized approximate inverses in the TT format}
\label{sec:ttfai}
\begin{figure}[tbp]
    \centering
    \includegraphics[width=0.75\textwidth]{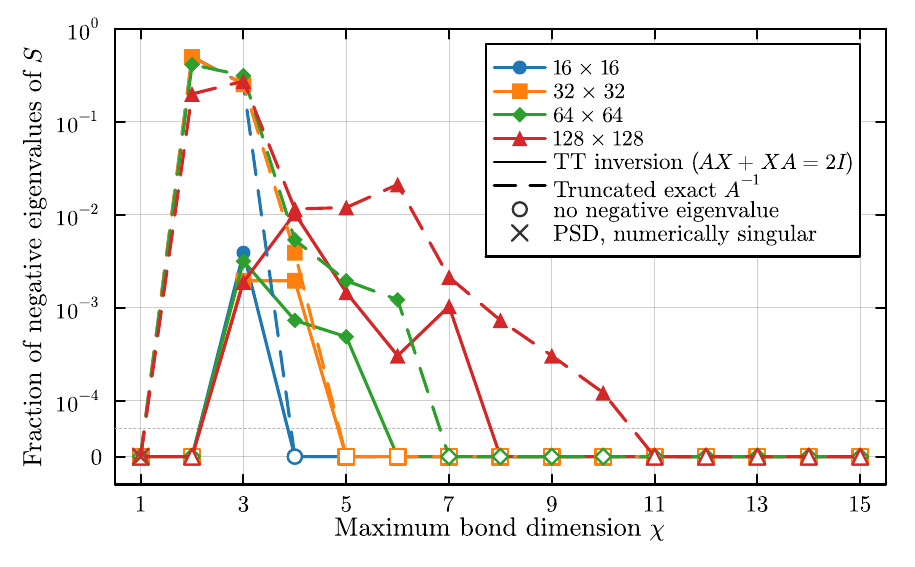}
    \caption{Fraction of negative eigenvalues of the symmetric part $S=(P+P^T)/2$ of two TT approximate inverses $P$
    of the pressure operator of Sec.~\ref{subsec:chorin_experiment}, as a function of their maximum bond dimension
    $\chi$, on grids of $16^2$, $32^2$, $64^2$ and $128^2$ points ($n = 256$, $1024$, $4096$ and $16384$). Solid lines: TT solution of
    $AX + XA = 2I$ \cite{oseledets2012}. Dashed lines: exact $A^{-1}$ converted to a TT operator and truncated to bond
    dimension $\chi$ by successive SVDs. Hollow markers on the bottom level indicate no negative eigenvalue, and crosses
    a positive semidefinite (PSD) but numerically singular $S$.}
    \label{fig:spd_nneg}
\end{figure}

Having reviewed the necessary background on the TT format and TT-based solvers, we now return to our main focus 
of accelerating the CG method. As mentioned, in order to precondition this solver, we seek an SPD approximate inverse to the operator 
of our linear system. Operator inversion techniques were already introduced in \cite{oseledets2010, oseledets2012} 
for the TT format, but even though one can obtain approximate inverses with great fidelity in terms of operator norm, 
the SPD property is not guaranteed by construction, possibly leading to breakdowns in the CG pipeline. 
To illustrate this, we consider the inverse of a two-dimensional Laplacian with
Dirichlet--Neumann boundary conditions in the $x$ direction and periodic boundary
conditions in the $y$ direction, represented in the QTT format of
Sec.~\ref{sec:tensor_trains}. Such a setting will be considered later in Sec.~\ref{subsec:chorin_experiment}. We approximate this inverse by TT operators
of increasing maximum bond dimension, by means of the TT inversion method of Ref.~\cite{oseledets2012} and by performing direct inversion and then converting 
the resulting matrix to a TT operator with a capped bond dimension. Fig.~\ref{fig:spd_nneg} shows the fraction of
negative eigenvalues of the symmetric part of each approximation. In both cases, truncation moves the smallest 
eigenvalues of the symmetric part of the inverse below zero. Truncating a TT operator controls the norm error but not 
the sign of its eigenvalues, and definiteness is only guaranteed once this error falls below 
$\lambda_{\min}(A^{-1}) = 1/\lambda_{\max}(A)$, which corresponds to a relative accuracy of order $1/\kappa(A)$.
This tolerance gets tighter as the grid gets refined, and, accordingly, the range of bond dimensions with 
an indefinite operator widens. On the grids of Sec.~\ref{subsec:chorin_experiment}, where dense spectra are out of reach, the 
local projections of the approximate inverses computed with the TT solver (Fig.~\ref{fig:projectedA}), whose spectra lie within the global one, 
are indefinite for $\chi$ up to 12 on the $256^2$ grid, and for every $\chi$ tested, up to 20, on the $1024^2$ one. 
Direct approximate inverses are therefore unreliable for PCG, precisely when aggressive truncation 
is needed to keep their cost low.

To enforce this property by construction in the TT format, we first turn to the theory of 
FSAIs, introduced by Kolotilina and 
Yeremin \cite{kolotilina1993}. The core idea behind FSAIs is to look for the inverse not of the operator itself, but of 
its Cholesky factor. Since the preconditioner is built in factorized form, it is positive semidefinite by construction, 
and definite whenever the factor is non-singular. 

We start by denoting the Cholesky decomposition of our operator $A$ and its approximate inverse $M^{-1}\approx A^{-1}$ as follows
\begin{gather}
    A = L_A L_A^T, \qquad M^{-1} = G_L^T G_L.
\end{gather}
With this description, $G_L$ is also an approximate inverse to 
$L_A$ and, in the exact case $G_L=L_A^{-1}$, it is lower triangular, like $L_A$. The latter statement can be shown by
applying column-wise forward substitution, for example. From this, we can frame the task of finding 
this $G_L$ factor as an optimization of the following cost function
\begin{align}
    \mathcal{L}(G_L) &= \|I - G_LL_A\|_F^2 = \text{Tr}\left\{(I - G_LL_A)(I - G_LL_A)^T\right\}\\\notag
    &=\text{Tr}\left\{I\right\} - \text{Tr}\left\{G_LL_A\right\} - \text{Tr}\left\{L_A^TG_L^T\right\} + \text{Tr}\left\{G_LL_AL_A^TG_L^T\right\},
\end{align}
whose stationary points with respect to the individual entries of $G_L$ satisfy
\begin{equation}
    \label{eq:fsai_linearsystem_withLa}
    (G_L A)_{ij} = (L_A^T)_{ij}\;, i \geq j.
\end{equation}
The condition $i \geq j$ reflects that the entries with respect to which we are differentiating 
only exist in the lower-triangular part. Since $L_A^T$ is upper triangular, the right-hand side vanishes for $i>j$, 
so that, up to the diagonal factor $\mathrm{diag}(L_A)$, $G_L$ coincides with the lower triangular approximate inverse 
of $A$ defined by $(\tilde G_L A)_{ij}=\delta_{ij}$, $i\ge j$, and the off-diagonal entries of $L_A$ are not needed \cite{kolotilina1993}.
This is a crucial property, derived from the triangular properties of these different matrices, that allows us to 
obtain $G_L$ without the need to compute the Cholesky factorization of $A$. To do this, we proceed as follows. We first 
compute a matrix $\tilde{G}_L$ that satisfies
\begin{equation}
\label{eq:fsai_linearsystem}
    (\tilde{G}_LA)_{ij} = \delta_{ij}\;, i \geq j.
\end{equation}
Then, we set $G_L :=D\tilde{G}_L$, with the diagonal matrix $D$ chosen such that the preconditioned matrix has a unit
diagonal, $(G_LAG_L^T)_{ii} = 1$ \cite{kolotilina1993}. Since $(G_LAG_L^T)_{ii} = D_{ii}^2(\tilde{G}_LA\tilde{G}_L^T)_{ii}$,
this gives
\begin{equation}
\label{eq:fsai_rescaling}
    D = \left[\text{diag}(\tilde{G}_L A \tilde{G}_L^T)\right]^{-\frac{1}{2}} = \left[\text{diag}(\tilde{G}_L)\right]^{-\frac{1}{2}},
\end{equation}
where this choice can be shown to minimize the Kaporin number of the preconditioned matrix $G_L A G_L^T$, which is a
specific performance measure for the CG method \cite{Kaporin1994}. The second equality follows from
$(\tilde{G}_LA\tilde{G}_L^T)_{ii} = \sum_{j\le i}(\tilde{G}_LA)_{ij}(\tilde{G}_L)_{ij} = (\tilde{G}_L)_{ii}$, where only
$j \le i$ contributes because $\tilde{G}_L$ is lower triangular, and where Eq.~\eqref{eq:fsai_linearsystem} sets
$(\tilde{G}_LA)_{ij} = \delta_{ij}$ for those entries. Eq.~\eqref{eq:fsai_rescaling} is what makes the rescaling cheap, since only the diagonal of
$\tilde{G}_L$ is needed, but it holds only if $\tilde{G}_L$ satisfies Eq.~\eqref{eq:fsai_linearsystem} exactly. We
return to this point in Lemma~\ref{lem:diag_positivity}, where $\tilde{G}_L$ is computed by an approximate solver.

\subsection{Obtaining \texorpdfstring{$\tilde{G}_L$}{G L} in the tensor train format}
\label{subsec:computing_G}
The solution of the linear system from Eq.~\eqref{eq:fsai_linearsystem}, when working with matrices, can be decomposed 
into independent linear systems, one for each row, that can be solved in a massively parallel way 
and that fully respect the lower triangular properties of the solution. On top of this, within these linear systems, 
specific sparsity constraints can be trivially imposed, keeping the cost of this factor under control at the expense 
of approximation quality \cite{janna2015}.

In the TT format, however, these same triangular properties 
cannot be imposed in such a straightforward way, since one typically works with individual pieces of a factorization instead of full matrices. In order to still end up with the desired construction, we will proceed 
by means of a projection step that globally modifies the linear system in order to restrict the possible solutions to
the subspace of lower-triangular operators.

\subsubsection{Construction of a triangularity mask in the TT format}

We assume that the operator $A$ is given as a TT operator 
in which each of the $N$ sites has a local dimension of $d_k$, where $k=1,\dots,N$, and introduce the following 
$d_k\times d_k$ blocks, given in terms of their entries:
\begin{equation}
\begin{aligned}
\mathbb{O}^{(d_k)}_{ij} &:= 1,\\
\mathbb{I}^{(d_k)}_{ij} &:= \delta_{ij},\\
\mathbb{L}^{(d_k)}_{ij} &:= \Theta(i-j), \\
\mathbb{J}^{(d_k)}_{ij} &:= \Theta(i-j-1) =\mathbb{L}^{(d_k)}_{ij} - \mathbb{I}^{(d_k)}_{ij}. \\
\end{aligned}
\qquad i,j=1,\dots,d_k,
\end{equation}
where $\Theta$ is the Heaviside step function, which is equal to one for non-negative arguments, and zero otherwise, 
and $i$ and $j$ denote the rows and columns, respectively. In simpler terms, this translates into $\mathbb{O}$ 
being a matrix full of ones, $\mathbb{I}$ being an identity and $\mathbb{L}$ and 
$\mathbb{J}$ being lower triangular matrices of ones, including and excluding the diagonal, respectively. We also introduce the matrix $P_L$ that 
serves as a mask for lower triangular matrices in the full space, given by
\begin{equation}
\left(P_L^{(1:N)}\right)_{ij} := \Theta(i-j) \qquad i,j=1,\dots, \prod_{k=1}^N{d_k},
\end{equation}
where the $1:N$ range is used to indicate precisely that its indices span dimensions from $1$ to $N$. If this mask
only applies to a single dimension $d_k$, we will simply denote it as $P_L^{(k)}$. The crucial property of this $P_L$ 
matrix is that it can be built from the $\mathbb{O}$, $\mathbb{I}$, $\mathbb{L}$ and $\mathbb{J}$ matrices as
a TT operator of bond dimension 2:
\begin{equation}
    P_L^{(1:N)} = \left[ \begin{matrix} \mathbb{I}^{(d_1)} & \mathbb{J}^{(d_1)}  \end{matrix} \right] \mathop{\bowtie}\limits_{k=2}^{N-1} \left(\left[ \begin{matrix} \mathbb{I}^{(d_k)} & \mathbb{J}^{(d_k)} \\  & \mathbb{O}^{(d_k)} \end{matrix} \right]\right) \bowtie \left[ \begin{matrix} \mathbb{L}^{(d_N)} \\ \mathbb{O}^{(d_N)} \end{matrix} \right].
\label{eq:lower_triangular_mask}
\end{equation}
The proof of this closely follows the approach taken in \cite{Kazeev2012}. We first consider the following recursive 
relationship, in which we relate $P_L^{(1:N)}$ and $P_L^{(2:N)}$ by adding a new dimension from the left:
\begin{align}
    P_L^{(1:N)} &= \left(\begin{array}{ccccc}
        P_L^{(2:N)} & & & & \\
        \mathbb{O}^{(d_2 \times \dots \times d_N)} & P_L^{(2:N)} & & & \\
        \vdots & \ddots & \ddots & & \\
        \mathbb{O}^{(d_2 \times \dots \times d_N)} & \dots & \mathbb{O}^{(d_2 \times \dots \times d_N)} & P_L^{(2:N)} & \\
        \mathbb{O}^{(d_2 \times \dots \times d_N)} & \dots & \dots & \mathbb{O}^{(d_2 \times \dots \times d_N)} & P_L^{(2:N)}
    \end{array}\right)\\
    & = \left[ \begin{matrix} \mathbb{I}^{(d_1)} & \mathbb{J}^{(d_1)}  \end{matrix} \right] \bowtie \left[ \begin{matrix} P_L^{(2:N)} \\ \mathbb{O}^{(d_2 \times \dots \times d_N)} \end{matrix} \right].
\label{eq:first projector decomposition}
\end{align}
Alongside this, we make use of the fact that 
$\mathbb{O}^{(d_1 \times \dots \times d_N)} = \mathbb{O}^{(d_1)} \otimes \mathbb{O}^{(d_2 \times \dots \times d_N)}$, 
which results in the following expansion of the identity from Eq.~\eqref{eq:first projector decomposition}:
\begin{equation}
\left[\begin{matrix} P_L^{(1:N)} \\ \mathbb{O}^{(d_1 \times \dots \times d_N)} \end{matrix}\right] = \left[ \begin{matrix} \mathbb{I}^{(d_1)} & \mathbb{J}^{(d_1)} \\ & \mathbb{O}^{(d_1)} \end{matrix} \right] \bowtie \left[ \begin{matrix} P_L^{(2:N)} \\ \mathbb{O}^{(d_2 \times \dots \times d_N)} \end{matrix} \right].
\end{equation}
Applying this identity recursively to the rightmost factor gives
\begin{align}
\left[\begin{matrix} P_L^{(1:N)} \\ \mathbb{O}^{(d_1 \times \dots \times d_N)} \end{matrix}\right] & = \left[ \begin{matrix} \mathbb{I}^{(d_1)} & \mathbb{J}^{(d_1)} \\ & \mathbb{O}^{(d_1)} \end{matrix} \right] \bowtie \left[ \begin{matrix} P_L^{(2:N)} \\ \mathbb{O}^{(d_2 \times \dots \times d_N)} \end{matrix} \right]\\\notag
& = \left[ \begin{matrix} \mathbb{I}^{(d_1)} & \mathbb{J}^{(d_1)} \\ & \mathbb{O}^{(d_1)} \end{matrix} \right] \bowtie \left[ \begin{matrix} \mathbb{I}^{(d_2)} & \mathbb{J}^{(d_2)} \\ & \mathbb{O}^{(d_2)} \end{matrix} \right] \bowtie \left[ \begin{matrix} P_L^{(3:N)} \\ \mathbb{O}^{(d_3 \times \dots \times d_N)} \end{matrix} \right]\\\notag
& = \dots = \mathop{\bowtie}\limits_{k=1}^{N-1} \left(\left[ \begin{matrix} \mathbb{I}^{(d_k)} & \mathbb{J}^{(d_k)} \\  & \mathbb{O}^{(d_k)} \end{matrix} \right]\right) \bowtie\left[ \begin{matrix} P_L^{(N)} \\ \mathbb{O}^{(d_N)} \end{matrix} \right]\notag
\end{align}
Identifying $P_L^{(N)} = \mathbb{L}^{(d_N)}$, the first block row is precisely Eq.~\eqref{eq:lower_triangular_mask}, completing the proof.

\subsubsection{Projection-based triangularity constraint}
\begin{figure}[tpb]
    \centering
    \begin{subfigure}[b]{0.45\textwidth}
        \centering
        \includegraphics[width=\textwidth]{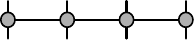}
        \caption{}
        \label{fig:mpo}
    \end{subfigure}\hfill
    \begin{subfigure}[b]{0.45\textwidth}
        \centering
        \includegraphics[width=\textwidth]{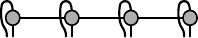}
        \caption{}
        \label{fig:vectorized_mpo}
    \end{subfigure}

    \vspace{0.5em}
    \begin{subfigure}[b]{0.45\textwidth}
        \centering
        \includegraphics[width=\textwidth]{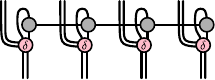}
        \caption{}
        \label{fig:pl_projector}
    \end{subfigure}
    \caption{Diagrammatic representations of the lower-triangular mask $P_L$ (a), its vectorized form
    (b), and the lower triangular mask projector $\mathcal{P}_L$ (c).}
    \label{fig:Projector processing pipeline}
\end{figure}

Now that we have a low-bond-dimension mask for lower triangular matrices, we can turn our attention back to the matrix 
linear system defined in Eq.~\eqref{eq:fsai_linearsystem_withLa}. First, we reframe it as a regular linear system by vectorizing 
the matrix space, as in the matrix inversion problem presented in Ref.~\cite{oseledets2012}
\begin{equation}
    (I\otimes A) \text{vec}(G_L) = \text{vec}(L_A^T).
    \label{eq:unprojected fai linear system}
\end{equation}
Here vec stacks the rows of its argument, so that $(I\otimes A)\,\text{vec}(G_L) = (I\otimes A^T)\,\text{vec}(G_L)= \text{vec}(G_L A)$ and $A$ 
acts on each row of $G_L$ separately, as in Eq.~\eqref{eq:fsai_linearsystem_withLa}.
As it is, this linear system has two problems. First, it involves the Cholesky factor of our operator $A$, which we cannot 
easily access in the TT format. Second, following the ALS-based pipeline of \cite{oseledets2012}, we would be solving the system as a whole, and, therefore,
we would not be imposing any lower-triangularity on the solution $G_L$. In order to tackle both issues, we 
can take the mask $P_L^{(1:N)}$ and obtain a projector $\mathcal{P}_L$ out of it in this vectorized space by first 
vectorizing it and then building a diagonal operator. This results in an operator that has a $1$ in its diagonal
if the entry belongs to the lower-triangular part and 0 otherwise. For illustrative purposes
the matrix version of this process for a mask of size $2 \times 2$ would look like
\begin{equation*}
    P_L = \left(\begin{array}{cc}
         1 & 0   \\
         1 & 1 
    \end{array}\right) \rightarrow \text{vec}(P_L) = \left[\begin{array}{c}
           1\\
           0\\
           1\\
           1\\
    \end{array}\right] \rightarrow \mathcal{P}_L = \left(\begin{array}{cccc}
         1 & 0 & 0 & 0  \\
         0 & 0 & 0 & 0 \\
         0 & 0 & 1 & 0 \\
         0 & 0 & 0 & 1 \\
    \end{array}\right) \,,
\end{equation*}
and Fig.~\ref{fig:Projector processing pipeline} contains the diagrammatic description of this process for $N=4$ sites in the 
TT format. Notice that the last step of this TT pipeline, which sets the vectorized $P_L^{(1:N)}$ in the diagonal 
of the expanded operator, involves simply contracting delta tensors site-wise, whose entries are equal to 1 when all the 
indices match in value and zero otherwise. The resulting projector, therefore, has the same low bond dimension as the 
starting lower-triangular mask.

Once we have the $\mathcal{P}_L$ object in TT form, we can project the linear system from Eq.~\eqref{eq:unprojected fai linear system} 
in order for the global solution to populate the lower triangular subspace:
\begin{equation}
    \label{eq:projected_fsai_linearsys}
    \left[\mathcal{P}_L(I\otimes A)\mathcal{P}_L + \mu(\mathcal{I} - \mathcal{P}_L)\right] \text{vec}(G_L) = \mathcal{P}_L\text{vec}(L_A^T) = \text{vec}(\text{diag}(L_A^T)),
\end{equation}
where $\mathcal{I}$ denotes an identity in the vectorized space, to distinguish it from the regular $I$ identity that has the same dimensions as $A$.
Let us review what each term in this new linear system does. The first one, involving the projector twice, 
sets the rows and columns of $I\otimes A$ that interact with the entries of the vectorized operator $G_L$ that 
live in the upper-triangular part to 0. The second term places a $\mu$ in the diagonal entries corresponding to those 
same upper-triangular entries. Finally, in the right-hand side, the projector, acting on the upper-triangular $L_A^T$ 
results in only its diagonal surviving. This construction, therefore, enforces the lower-triangular property on $G_L$ 
variationally through the upper-triangular entries appearing in the linear system as $\mu G_{L_{(i<j)}}=0$. 
This is, of course, less robust than the sparse matrix approach that guarantees the triangularity property by construction, since
an approximate solver can still set those entries to non-zero values. Still, the important property that we get from 
this projected linear system is that in the right-hand side we recover the diagonal operator that allows us to proceed in two steps and avoid computing the 
Cholesky factor $L_A$, as covered in the beginning of this section.

We can therefore begin by solving the following linear system 
\begin{equation}
    \label{eq:projected_fsai_linearsys_withI}
    \left[\mathcal{P}_L(I\otimes A)\mathcal{P}_L + \mu(\mathcal{I} - \mathcal{P}_L)\right] \text{vec}(\tilde{G}_L) = \text{vec}(I) \,,
\end{equation}
with the TT solver of choice, and then rescale the TT operator $\tilde{G}_L$ with another one $D$ given by 
the element-wise operation
\begin{equation}
    \label{eq:rescaling D}
    D = \left[\text{diag}(\tilde{G}_L)\right]^{-\frac{1}{2}} \,.
\end{equation}
There are two important things worth noting in this pipeline, however. First, even though the value of $\mu$ does not affect the exact solution of the linear system of Eq.~\eqref{eq:projected_fsai_linearsys_withI}, as discussed in \ref{app:row-wise splitting}, it does affect the condition number of the effective operator, given by
\begin{equation}
    \kappa\left(\mathcal{P}_L(I\otimes A)\mathcal{P}_L + \mu(\mathcal{I} - \mathcal{P}_L)\right) = \frac{\max(\mu, \lambda_{\max}(A))}{\min(\mu, \lambda_{\min}(A))}.
\end{equation}
A $\mu$ value lying on the upper part of the spectrum of $A$, therefore, will provide a large penalization to the upper-triangular entries without affecting the performance of the CG method when solving the linear system. $\lambda_{\max}(A)$ can be estimated analytically for some particular cases or by means of a few iterations of the Krylov method of choice, for example. 
Second, regarding the rescaling step, Eq.~\eqref{eq:rescaling D}  is the shortcut in Eq.~\eqref{eq:fsai_rescaling}, in which the energy-based scaling $(\tilde{G}_LA\tilde{G}_L^T)_{ii}$ is replaced by the diagonal entry $(\tilde{G}_L)_{ii}$.
In classical FSAI each row of $\tilde{G}_L$ is obtained by solving a small and dense SPD system, typically with a direct
method. The replacement is then exact, and the diagonal entries are positive for any sparsity pattern that contains the
main diagonal, since they are diagonal entries of inverses of SPD principal submatrices of $A$
\cite[Lemma 2.1]{kolotilina1993}. Here, in contrast, Eq.~\eqref{eq:projected_fsai_linearsys_withI} is solved as a
whole by an approximate TT solver, so neither property is inherited automatically. In the sparse setting, the same
situation arises whenever the rows are not exact solutions of their restricted systems, and the shortcut is then 
avoided by computing the scaling from the energy of each row, $D_{ii} = (\tilde{\mathbf{g}}_i^TA\tilde{\mathbf{g}}_i)^{-1/2}$, 
which is positive for any non-zero row \cite{janna2015}. In the TT format, extracting
$\text{diag}(\tilde{G}_LA\tilde{G}_L^T)$ requires contracting three tensor
train operators, which is an expensive operation that can require a large bond dimension,
whereas $\text{diag}(\tilde{G}_L)$ is
obtained with a single contraction of delta tensors (see Fig.~\ref{fig:diagG}). We therefore keep the cheaper
Eq.~\eqref{eq:rescaling D} (see, however, the remark at the end of Sec.~\ref{subsec:computing_D}), and quantify
instead how far its diagonal entries can drift from the exact ones, and when they can change sign.

To state when the diagonal entries of $\tilde{G}_L$ stay positive, in order for Eq.~\eqref{eq:rescaling D} to be well defined, we denote by $\mathcal{T} = \mathcal{P}_L(I\otimes A)\mathcal{P}_L + \mu(\mathcal{I} -
\mathcal{P}_L)$ the operator of Eq.~\eqref{eq:projected_fsai_linearsys_withI}, by $\tilde{G}_L^\star$ its exact
solution, and by $\tilde{G}_L$ any approximation to it, such as the output of a TT solver with a finite bond dimension and
tolerance. Since $\mathcal{P}_L$ acts on each row of $\tilde{G}_L$ separately, $\mathcal{T}$ is block diagonal, with one
block $\mathcal{T}_i$ per row (see \ref{app:diag_positivity}). We write $\mathbf{x}_i$ and $\mathbf{x}_i^\star$ for the
$i$-th rows of $\tilde{G}_L$ and $\tilde{G}_L^\star$, $g_i^\star = (\tilde{G}_L^\star)_{ii}$,
$\|\mathbf{v}\|_{\mathcal{T}_i}^2 = \mathbf{v}^T\mathcal{T}_i\mathbf{v}$ for the energy norm of row $i$, and
\begin{equation}
    \label{eq:row energy error}
    \nu_i = \frac{\|\mathbf{x}_i - \mathbf{x}_i^\star\|_{\mathcal{T}_i}}{\|\mathbf{x}_i^\star\|_{\mathcal{T}_i}}
\end{equation}
for the relative energy-norm error of row $i$.

\begin{lemma}[Positivity of the diagonal of $\tilde{G}_L$]
\label{lem:diag_positivity}
Let $A$ be SPD with Cholesky factor $L_A$, and let $\mu > 0$. Then, for every row $i$:
\begin{enumerate}
    \item[(a)] The exact diagonal is positive: $g_i^\star = 1/(L_A)_{ii}^2$, and
    $\dfrac{1}{A_{ii}} \le g_i^\star \le \dfrac{1}{\lambda_{\min}(A)}$.
    \item[(b)] The computed diagonal is stable:
    $\left|\dfrac{(\tilde{G}_L)_{ii} - g_i^\star}{g_i^\star}\right| \le \nu_i$.
    \item[(c)] A sign change requires losing the whole row: if $(\tilde{G}_L)_{ii} \le 0$, then $\nu_i \ge 1$.
\end{enumerate}
\end{lemma}

\begin{corollary}
\label{cor:diag_fraction}
Let $\varepsilon_E = \|\mathrm{vec}(\tilde{G}_L - \tilde{G}_L^\star)\|_{\mathcal{T}} /
\|\mathrm{vec}(\tilde{G}_L^\star)\|_{\mathcal{T}}$ be the global relative energy-norm error. Then the rows with a
non-positive diagonal entry satisfy
\begin{equation}
    \label{eq:diag fraction}
    \frac{\sum_{i\,:\,(\tilde{G}_L)_{ii}\le 0}\, g_i^\star}{\sum_i g_i^\star} \le \varepsilon_E^2 .
\end{equation}
\end{corollary}
The proofs are given in \ref{app:diag_positivity}, together with a discussion of why sign changes are not observed in
practice and of what the lemma does not guarantee. 

It is worth noting, however, that even if all the entries of $\text{diag}(\tilde{G}_L)$ are positive, the subsequent 
steps to obtain the rescaled $G_L$ also introduce errors. The calculation of the element-wise operation from 
Eq.~\eqref{eq:rescaling D} is approximate in itself, as discussed in the following section, and truncation errors 
are introduced during the contraction of the TT operators representing $D$ and $\tilde{G}_L$. Importantly, 
however, none of these errors and approximations can make the preconditioner indefinite, as long as the 
sequential application of $G_L$ and its transpose is performed in an exact manner. 

By setting a maximum allowed bond dimension for $\tilde{G}_L$, $D$ and the 
product of the two, then, this pipeline produces an approximately lower-triangular, low-rank TT approximation to $L_A^{-1}$.
The bond dimension of the outcome $G_L$ thus acts as an analog to the sparsity pattern one can
impose when working with sparse matrices in the classical setting. Regardless of whether $L_A^{-1}$ admits an accurate low-rank
representation, it gives us the ability to externally control the cost of storing and applying our preconditioner at the 
expense of preconditioning performance. 

\subsection{Computing the rescaling factor \texorpdfstring{$D$}{D} in the TT format}
\label{subsec:computing_D}
\begin{figure}[t]
    \centering
    \begin{subfigure}[b]{0.45\textwidth}
        \centering
        \includegraphics[width=\textwidth]{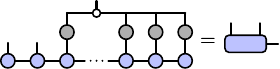}
        \caption{Initial sketch.}
        \label{fig:first_sketch}
    \end{subfigure}\hfill
    \begin{subfigure}[b]{0.45\textwidth}
        \centering
        \includegraphics[width=\textwidth]{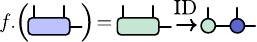}
        \caption{Element-wise function application and interpolative decomposition (ID).}
        \label{fig:id}
    \end{subfigure}
        \begin{subfigure}[b]{0.45\textwidth}
        \centering
        \includegraphics[width=\textwidth]{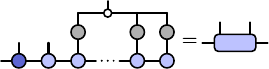}
        \caption{Second sketch, from reconstructed TT.}
        \label{fig:second_sketch}
    \end{subfigure}
    \caption{Core steps of the RSI pipeline for implementing element-wise functions of TTs. Step (a)
    performs a random sketch of the TT by means of a Khatri--Rao product of random matrices, depicted by the gray nodes. 
    Then, in step (b), the desired function is applied element-wise (notice the \textit{dot} after $f$, which represents the Julia language convention for 
    element-wise operations) on the sketched array and an interpolative 
    decomposition is applied.  This recovers the TT format on the output and lets us reconstruct the original TT in the 
    basis defined by the output, denoted by the darker shade of blue. This allows us to repeat the process along the 
    chain, as shown in step (c), until the whole TT has been processed.}
    \label{fig:RSI}
\end{figure}
In order to compute the rescaling factor $D$ as given by Eq.~\eqref{eq:rescaling D}, we 
employ the Recursive Sketched Interpolation (RSI) pipeline introduced in \cite{meng2026}, which was presented mainly for performing 
element-wise multiplication of TTs, but that can also accommodate other kinds of operations, such as the square 
root inverse that appears in Eq.~\eqref{eq:rescaling D}. The core of this algorithm consists of mapping the TT to be processed onto just a subset of its sites or nodes by means of a random map. This process, known as 
\textit{sketching}, then allows one to perform the element-wise operation on an array of manageable size. After this operation
has been performed, the factorized format can be recovered by means of an interpolative decomposition, which also allows for the reconstruction of the original TT with matching indices, so that the operation can be repeated on the next subset of sites. By performing this operation throughout the network, then, one obtains the full TT that represents 
the processed data. Fig.~\ref{fig:RSI} summarizes the pipeline in diagrammatic terms. In the original formulation of RSI, however, these random sketches are sampled from a normal distribution, 
but the function that we want to apply only has support on positive entries. Assuming that the solver that tackled
Eq.~\eqref{eq:projected_fsai_linearsys_withI} converged to a meaningful result, the operator $\tilde{G}_L$ should
contain only positive numbers in its diagonal, as discussed in the previous section, but the contraction of Gaussian sketches would result in an array that is not 
necessarily composed of positive numbers only. In order to overcome this, we have replaced the Gaussian sketches with 
uniformly distributed sketches, which result in a sketched array that is strictly 
positive, allowing us to apply the function $f(x) = x^{-1/2}$ element-wise without issues. 

\begin{figure}[htbp]
    \centering
    \begin{subfigure}[b]{0.45\textwidth}
        \centering
        \includegraphics[width=\textwidth]{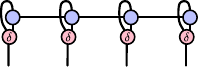}
        \caption{$\text{diag}(\tilde{G}_L)$ in TT format.}
        \label{fig:diagG}
    \end{subfigure}\hfill
    \begin{subfigure}[b]{0.45\textwidth}
        \centering
        \includegraphics[width=\textwidth]{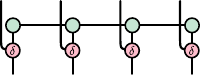}
        \caption{Rescaling $D$ in TT format.}
        \label{fig:Dmpo}
    \end{subfigure}
    \caption{Contraction of delta tensors on $\tilde{G}_L$ to obtain $\text{diag}(\tilde{G}_L)$ (a) in order to apply the function $f(x) = x^{-1/2}$
    element-wise via the RSI pipeline, and recovery of the TT operator form from the output of RSI (b).}
    \label{fig:Dpipeline}
\end{figure}

The processing pipeline to obtain $D$ is, therefore, the following. First, a TT of $\text{diag}(\tilde{G}_L)$ is 
obtained by contracting deltas onto the TT operator $\tilde{G}_L$, as shown in Fig.~\ref{fig:diagG}. This object is then passed to RSI 
in order to apply the inverse of the square root element-wise. The output TT is then contracted with deltas, 
this time to recover the operator form, as shown in Fig.~\ref{fig:Dmpo}. Once we have this rescaler TT operator, 
therefore, we can contract it onto $\tilde{G}_L$ and obtain the factor $G_L$ of the approximate inverse.

\begin{figure}[htbp]
    \centering
    \includegraphics[width=0.4\textwidth]{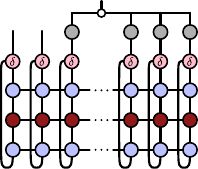}
    \caption{First RSI sketch of $\text{vec}(\text{diag}(\tilde{G}_LA\tilde{G}_L^T))$ kept as an uncontracted composite
    network. The cores of $\tilde{G}_L$ and $\tilde{G}_L^T$ (blue) and of $A$ (dark red) are joined site by site by delta
    tensors, and the random sketching matrices (gray) are contracted onto the resulting network as in Fig.~\ref{fig:first_sketch}, without ever forming the product $\tilde{G}_LA\tilde{G}_L^T$ as a TT operator.}
    \label{fig:energy_sketch}
\end{figure}

We note, however, that the cost argument given in Sec.~\ref{subsec:computing_G} against the energy-based scaling
applies mainly to forming $\tilde{G}_LA\tilde{G}_L^T$ explicitly. RSI only accesses its input through contractions with the
sketching matrices, so the input need not be available as a single, rounded TT. By keeping the object
$\text{vec}(\text{diag}(\tilde{G}_LA\tilde{G}_L^T))$ as an uncontracted composite network, the sketched matrices involved
in RSI to obtain $D$ could be obtained by contracting into it the random matrices without ever forming an explicit
TT operator, as shown in Fig.~\ref{fig:energy_sketch}, similarly to the evaluation of an expectation value, for example. 
The energy-based scaling could therefore be obtained without any error beyond that of RSI
itself. Denoting by $\chi_{G_L}$ and $\chi_A$ the bond dimensions of $\tilde{G}_L$ and $A$, the leading cost of building each two-site block of
Fig.~\ref{fig:second_sketch} would then grow by a factor of about $\chi_A\chi_{G_L}$. Building the MPO of $\tilde{G}_LA\tilde{G}_L^T$
explicitly and without truncation would instead cost $\mathcal{O}(d^3\chi_{G_L}^4\chi_A^2)$ per site, and its bond
dimension of $\chi_{G_L}^2\chi_A$ would raise the cost of each block by a factor of about $\chi_A^2\chi_{G_L}^2$ with respect to the direct rescaling of Eq.~\eqref{eq:rescaling D}. Moreover, since $A$ is SPD, these entries are positive for every non-zero row, so this energy-based variant would make $D$
well defined regardless of the accuracy of $\tilde{G}_L$. We nevertheless leave the exploration of this alternative for
future work. On the uniform grids with constant coefficients considered here, the exact diagonal of
$\tilde{G}_L$ is nearly constant, and a non-positive entry would require an error of the size of a whole row
concentrated on that row, which TT solvers have no reason to produce (see \ref{app:sign_changes}). The energy-based
variant would become relevant on strongly non-uniform grids or with strongly varying coefficients, for which this
argument no longer holds (see \ref{app:diag_limits}).

\section{Factorized preconditioning of DMRG-like algorithms}
\label{sec:Preconditioning ALS}
\begin{figure}[tpb]
    \centering
    \includegraphics[width=0.4\textwidth]{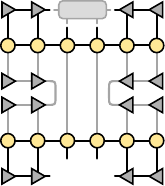}
    \caption{Projected factorized preconditioner $G_{L_{\text{loc}}}^T G_{L_{\text{loc}}} = U^T G_L^T U\, U^T G_L U$.}
    \label{fig:projectedGlTGl}
\end{figure}

Now that we have an SPD approximate inverse to our operator, we can move on to using it to improve the convergence of TT 
linear solvers based on the CG method. In order to use the factorized preconditioner within DMRG-like algorithms, however,
there is still some freedom in how the global preconditioner gets projected onto the local subspace. One could
contract the two factors into a single TT operator and proceed with the alternating optimization as covered in 
Sec.~\ref{subsec:tt_linear_solvers}. That, however, would substantially increase the bond dimension of the merged object, making 
the use of this preconditioner very costly, and truncating this product would break the positive-definite guarantees
of the factorized preconditioner, in a similar fashion to what is shown in Fig.~\ref{fig:spd_nneg}. 
The alternative approach we propose avoids these two issues by leveraging the fact that the two operators are not 
independent, but transposes of each other. By keeping track of just the projected factor $G_L$, 
we can also apply $G_L^T$ locally by transposing its components. In this sense, 
therefore, the cost of storing and updating the operators is the same as for a regular TT operator and, as the two factors of bond dimension $\chi_{G_L}$ are applied on top of $A$, a preconditioned local iteration costs about $1 + 2\chi_{G_L}/\chi_{A}$ unpreconditioned ones. For illustrative purposes, Fig.~\ref{fig:projectedGlTGl} showcases the contraction of the projected factorized preconditioner, composed of the 
projected $G_L$ and its transpose. 
It is worth noting, however, that proceeding in this way, in which one projects the 
factors and not the full preconditioner $G_L^T G_L$, introduces an intermediate $U U^T$ that leads to
a different preconditioner. 
In our experiments, we have not observed any qualitative differences in terms of 
iteration counts between this approach and the much more costly $U^T G_L^T G_L U$, meaning that this projection 
still yields a preconditioner of comparable quality that, more importantly, is also locally positive semidefinite, and definite whenever $U^T G_L U$ is non-singular.

\section{Numerical experiments}
\label{sec:experiments}

In this section we present the performance of TT-FAI preconditioners by analyzing their behavior in two
different elliptic problems: a high-dimensional anisotropic diffusion equation and the pressure Poisson equation of a
projection method for incompressible flows. The choice of second-order elliptic operators follows from three main
arguments. First, their discretizations are SPD. Second, their condition number grows polynomially as
$\kappa \sim \mathcal{O}(h^{-2})$, where $h$ is the grid spacing of the finite difference scheme of choice. Third,
they are of great significance in general science and engineering, appearing in a variety of problems in
electromagnetics, fluid and solid mechanics or quantitative finance.

For each of the two problems below, we benchmark the TT-FAI preconditioner on a DMRG-like MALS solver. In particular, the variant with residual-based truncation introduced in Ref.~\cite{oseledets2012}
has been used, as it proved more robust in meeting the prescribed residual bounds.
The TT-FAI preconditioner itself is, in every case, also constructed by
solving the projected linear system of Eq.~\eqref{eq:projected_fsai_linearsys_withI} with the regular MALS solver of Ref.~\cite{holtz2012}, as this
offered the best performance in preliminary tests, given the limited bond dimension budget allowed for $\tilde{G}_L$. The TT-FAI preconditioner is
compared, in every case, against an unpreconditioned baseline and, where applicable, against the preconditioner of Ref.~\cite{khoromskij2009b}. All solvers and preconditioners in this section are run with DifferenTT.jl \cite{DifferenTT}, and the scripts and data behind every figure and table are described in the Code and data availability statement. All runs used single MareNostrum 5 NG nodes (Nvidia Grace CPU superchip), with 20 BLAS threads and one Julia thread.

As for the choice of $\mu$, following the argument presented in Sec.~\ref{subsec:computing_G}, instead of estimating $\lambda_{\max}(A)$ numerically, in both experiments we set
it to the bound $\mu = \|A\|_\infty \ge \lambda_{\max}(A)$, which follows from Gershgorin's theorem \cite[App.~C.8]{Leveque2007}. For the finite difference operators on regular grids used in our benchmarks, this bound can be read directly from the stencil coefficients. For a central stencil
$\sum_k c_k u_{i+k}/h^2$ along each axis, every row sum of $|A|$ is bounded by $\sigma/h^2$ with
$\sigma = \sum_k |c_k|$, a bound that the boundary closures used here do not exceed. The resulting values are given
in each subsection below.

\subsection{High-dimensional anisotropic diffusion}
\label{subsec:anisotropy_experiment}

In this first experiment we evaluate the performance of the TT-FAI preconditioner in the non-quantized TT format against
the approximate inverse preconditioner introduced in Ref.~\cite{khoromskij2009b}. This preconditioner targets operators
in Kronecker-sum form, $\mathcal{L}_0 = \sum_{\ell=1}^N \mathcal{A}_\ell$, where each $\mathcal{A}_\ell$ is the SPD
matrix of a one-dimensional elliptic operator acting on the $\ell$-th dimension, and is given by the following expression:
\begin{gather}
    \label{eq:sinc preconditioners}
    B_M := \frac{1}{\lambda_0}\sum_{k=-M}^{M} c_k \bigotimes_{\ell=1}^N \exp\!\left(-\frac{t_k}{\lambda_0} \mathcal{A}_\ell\right) \approx \mathcal{L}_0^{-1},\\\notag
    t_k = e^{k\eta}, \quad c_k = \eta\, t_k, \quad \eta = \frac{\pi}{\sqrt{M}}, \quad
    \lambda_0 := \lambda_{\min}(\mathcal{L}_0) = \sum_{\ell=1}^N \lambda_{\min}(\mathcal{A}_\ell).
\end{gather}
The weights and nodes $(c_k, t_k)$ come from the Sinc quadrature of $1/x = \int_{\mathbb{R}} e^{s} e^{-e^{s}x}\,ds$
with step $\eta$, which is accurate for $x \in [1, R]$. For this reason the quadrature is applied to the rescaled
operator $\mathcal{L}_0/\lambda_0$, whose spectrum lies in $[1, \kappa(\mathcal{L}_0)]$ regardless of the scaling of
the $\mathcal{A}_\ell$ (e.g. their $1/h^2$ factor), and the result is scaled back by $1/\lambda_0$. Since the
eigenvalues of a Kronecker sum are the sums of the eigenvalues of its terms, $\lambda_0$ only requires the smallest
eigenvalue of each one-dimensional $\mathcal{A}_\ell$.
Each term of this sum is a rank-one TT operator, so $B_M$ has a bond dimension of at most $2M+1$, and since $c_k > 0$
and all the factors are SPD, so is $B_M$. Moreover, it can be shown that this approximate inverse converges to
$\mathcal{L}_0^{-1}$ exponentially fast in $\sqrt{M}$ \cite{khoromskij2009b}.

The $N$-dimensional discrete Laplacian is precisely of this form, and
therefore constitutes the ideal setting for Eq.~\eqref{eq:sinc preconditioners}, in which a handful of quadrature
terms already provide an almost exact inverse. On top of this, since the eigenvalues of a Kronecker sum are the sums of
the eigenvalues of its terms, with the same grid along every dimension, $\kappa(-\Delta_h) = \kappa(\mathcal{A}_\ell)$
does not depend on $N$, so increasing the number of dimensions does not make the Poisson equation any harder to
precondition. In order to benchmark the TT-FAI in a setting that departs from this ideal case in a controlled manner, we
instead consider the anisotropic diffusion equation
\begin{equation}
    \label{eq:anisotropic diffusion}
    -\nabla\cdot(K\nabla u) = f \ \text{ in } (0,1)^N, \quad u = 0 \ \text{ on } \partial(0,1)^N, \quad
    K = (1-\rho)\,\mathbb{I}^{(N)} + \rho\, \mathbf{1}\mathbf{1}^T,
\end{equation}
where $\mathbf{1}$ is the all-ones vector of length $N$, so that $\mathbf{1}\mathbf{1}^T$ is the
$N \times N$ matrix of ones. The resulting equicorrelated diffusion tensor $K$ has eigenvalue $1+(N-1)\rho$ along the diagonal direction
$\mathbf{1}/\sqrt{N}$ and $1-\rho$, with multiplicity $N-1$, on its orthogonal complement, and is SPD for
$-1/(N-1) < \rho < 1$. For $\rho > 0$, this describes a rotated anisotropy, with a strong diffusion along the main
diagonal of the hypercube and a weak one transverse to it, which is not aligned with any of the grid axes. In two
dimensions, this is the classic rotated anisotropy test case of the multigrid literature \cite{trottenberg2001}. In
higher dimensions, the same diffusion operator appears in problems driven by correlated Brownian motions \cite{oksendal2003}.

Since $\mathrm{diag}(K) = \mathbb{I}^{(N)}$, discretizing Eq.~\eqref{eq:anisotropic diffusion} with centered finite differences yields
\begin{equation}
    \label{eq:anisotropic operator}
    \mathcal{L}(\rho) = \sum_{\ell=1}^N \mathcal{A}_\ell - 2\rho \sum_{\ell<m} \mathcal{C}_\ell\, \mathcal{C}_m,
\end{equation}
where $\mathcal{C}_\ell$ is the centered first difference along the $\ell$-th dimension.
Even though the second term couples every pair of dimensions, $\mathcal{L}(\rho)$ admits an
exact MPO representation with a bond dimension of 3 for any $N$, given by
\begin{equation}
    \label{eq:anisotropic mpo}
    \mathcal{L}(\rho) = \left[ \begin{matrix} \mathbb{I}^{(d_1)} & \mathcal{C}^{(d_1)} & \mathcal{A}^{(d_1)} \end{matrix} \right]
    \mathop{\bowtie}\limits_{\ell=2}^{N-1} \left(\left[ \begin{matrix} \mathbb{I}^{(d_\ell)} & \mathcal{C}^{(d_\ell)} & \mathcal{A}^{(d_\ell)} \\ & \mathbb{I}^{(d_\ell)} & -2\rho\, \mathcal{C}^{(d_\ell)} \\ & & \mathbb{I}^{(d_\ell)} \end{matrix} \right]\right)
    \bowtie \left[ \begin{matrix} \mathcal{A}^{(d_N)} \\ -2\rho\, \mathcal{C}^{(d_N)} \\ \mathbb{I}^{(d_N)} \end{matrix} \right],
\end{equation}
whose entries are the $d_\ell\times d_\ell$ one-dimensional identity $\mathbb{I}^{(d_\ell)}$, first difference $\mathcal{C}^{(d_\ell)}$ and
second difference $\mathcal{A}^{(d_\ell)}$ matrices acting on each site. This is easily verified, for example, by means of 
finite automata \cite{Crosswhite2008}. For $\rho = 0$, it reduces to the bond dimension 2 MPO of the Laplacian 
\cite{Kazeev2012}, $\mathcal{L}(0) = -\Delta_h$.

The key property of this problem is that its distance from the separable regime of Eq.~\eqref{eq:sinc preconditioners}
can be tuned through $\rho$. Writing $\mathcal{L}(\rho) = \sum_\ell (\mathcal{A}_\ell - \mathcal{C}_\ell^T\mathcal{C}_\ell)
+ \sum_{\ell,m} K_{\ell m}\, \mathcal{C}_\ell^T \mathcal{C}_m$ and using that $0 \preceq \mathcal{C}^T\mathcal{C}
\preceq \mathcal{A}$ for the one-dimensional matrices, one finds (see \ref{app:anisotropy_floor} for a step-by-step
derivation) that $\lambda_{\min}(K)\,\mathcal{L}_0 \preceq
\mathcal{L}(\rho) \preceq \lambda_{\max}(K)\,\mathcal{L}_0$, where $\mathcal{L}_0 = -\Delta_h$ is the separable
operator obtained by dropping the off-diagonal part of $K$. Here, $X \preceq Y$ indicates that $\mathbf{v}^TX\mathbf{v} \leq \mathbf{v}^TY\mathbf{v},\,\forall\mathbf{v}$. Hence, for $\rho \ge 0$,
\begin{equation}
    \label{eq:anisotropy floor}
    \kappa(\mathcal{L}_0^{-1}\mathcal{L}) \le \frac{1+(N-1)\rho}{1-\rho} =: F.
\end{equation}
The quantity $\kappa(\mathcal{L}_0^{-1}\mathcal{L})$ acts as a floor for any preconditioner that approximates $\mathcal{L}_0^{-1}$. The TT-FAI, on the other hand, is built directly from
$\mathcal{L}(\rho)$, and is therefore not bound by this floor.

In order to compare different numbers of dimensions at an equal difficulty for separable preconditioners, rather than
sweeping $\rho$ directly, we fix the bound $F$ of Eq.~\eqref{eq:anisotropy floor} and invert it for $\rho$
\begin{equation}
    \label{eq:rho from floor}
    \rho(F, N) = \frac{F-1}{F+N-1},
\end{equation}
with $F \in \{1, 2, 4, 10, 30, 100\}$. The case $F=1$ ($\rho = 0$) recovers the Poisson equation, the most favorable
setting for $B_M$, while increasing $F$ progressively moves the problem away from it. As $F$ is only an upper bound, we also compute the discrete floor $\kappa(\mathcal{L}_0^{-1}\mathcal{L})$ on the benchmark grid of 15 points per
dimension for $N = 5$, in full format, from the Lanczos tridiagonal matrix associated with an exactly
$\mathcal{L}_0$-preconditioned CG run \cite{saad2003}, as this is not feasible for $N=10$. The analytic bound is not
attained on this coarse grid, since the modes required to reach the extreme eigenvalues of $K$ are not fully resolved,
but the discrete floor still grows monotonically, from 1 for $F=1$ to 20.4 for $F=100$.

As the right-hand side, we choose $f \equiv 1$, a rank-one source that, unlike the manufactured solutions commonly used
in TT benchmarks, does not add any structure of its own to the problem, so all the complexity of the solution comes
from the operator. Indeed, on the same grid with $N=5$, the maximum TT rank of the exact solution at a relative
accuracy of $10^{-4}$ grows steadily with the coupling, from 4 for $F=1$ to 27 for $F=100$, so the solves become
increasingly demanding as $F$ grows.

We compare three kinds of preconditioners against the unpreconditioned baseline. First, the preconditioner
$B_M$ of Eq.~\eqref{eq:sinc preconditioners} built on $\mathcal{L}_0$, which is the same for every $F$. Second, the TT-FAI
built on $\mathcal{L}(\rho)$. And finally, as a control, the TT-FAI built on the separable $\mathcal{L}_0$ but used to
precondition $\mathcal{L}(\rho)$, which is subject to the same floor as $B_M$. Comparing the two variants,
therefore, isolates whether the performance of the TT-FAI stems from capturing the coupling between dimensions or just
from approximating $\mathcal{L}_0^{-1}$. 
From the analysis above, $B_M$ is expected to be the best choice for $F=1$, where it is close to exact at a negligible
construction cost. As $F$ grows, its iteration count should saturate with $M$ at a level set by the floor, while the
TT-FAI of $\mathcal{L}(\rho)$ can, in principle, go below it at the price of a larger bond dimension. 

For the Gershgorin value of $\mu$ involved in the TT-FAI of $\mathcal{L}(\rho)$, the second differences of $\mathcal{L}$ contribute $4/h^2$ per dimension and each of the $N(N-1)/2$ mixed
terms $2\rho\,\mathcal{C}_\ell\mathcal{C}_m$ contributes $2\rho/h^2$, so that
\begin{equation}
\label{eq:anisotropy_gershgorin}
    \mu = \|\mathcal{L}(\rho)\|_\infty = \frac{N\left(4 + \rho(N-1)\right)}{h^2},
\end{equation}
which reduces to $4N/h^2$ for the TT-FAI of $\mathcal{L}_0$. 

\begin{table}[htbp]
    \centering
    \begin{tabular}{@{}ccccc@{}}
        \toprule
        $N$ & Max. bond dim. & Max. sweeps & $\mu$, TT-FAI of $\mathcal{L}$ & $\mu$, TT-FAI of $\mathcal{L}_0$ \\
        \midrule
        10 & 150 & 30 & $7.8\times10^3$--$2.4\times10^4$ & $7.8\times10^3$ \\
        20 & 300 & 50 & $1.6\times10^4$--$7.8\times10^4$ & $1.6\times10^4$ \\
        \bottomrule
    \end{tabular}
    \caption{MALS parameters for the anisotropic diffusion benchmark, with $F \in \{1,2,4,10,30,100\}$ and
    $f \equiv 1$, with 15 grid points per dimension. The outer relative residual tolerance is $10^{-5}$, with a residual-based truncation tolerance of $5 \times 10^{-6}$ and a cutoff of
    $10^{-13}$. The local CG solver runs up to a relative tolerance of $10^{-6}$ with at most 3000 iterations per
    local problem. The $B_M$ preconditioner uses $M \in \{3,4,5,7,10\}$, the TT-FAI of $\mathcal{L}(\rho)$
    $\chi \in \{4,6,8,10\}$ and the TT-FAI of $\mathcal{L}_0$ $\chi = 10$. Both TT-FAI variants use
    the $\mu$ value given by Eq.~\eqref{eq:anisotropy_gershgorin} and are built with at most 10 MALS sweeps at local
    tolerance $10^{-5}$, 2000 maximum CG iterations and a global residual tolerance of $10^{-1}$.}
    \label{tab:anisotropy_params}
\end{table}

Table~\ref{tab:anisotropy_params} lists the parameters of the MALS benchmark, for which we consider
$N \in \{10, 20\}$. Every run reaches the outer tolerance in exactly two sweeps, and every TT-FAI build reaches its
residual tolerance after a single one. Tighter build tolerances did not improve the results in preliminary tests. Every configuration has been repeated 10 times, each with a different random
initial guess and, for the TT-FAI, a new random sketch for $D$, and the figures below report the
median and the interquartile range over these repeats.

\begin{figure}[htbp]
    \centering
    \begin{subfigure}[b]{0.49\textwidth}
        \centering
        \includegraphics[width=\textwidth]{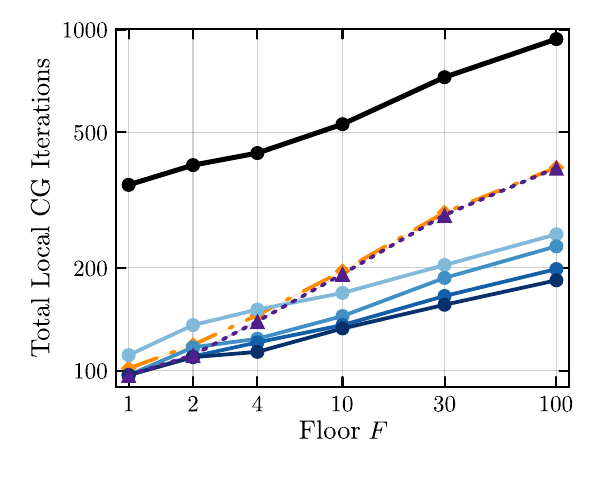}
        \caption{$N=10$}
        \label{fig:mals_iters_N10}
    \end{subfigure}\hfill
    \begin{subfigure}[b]{0.49\textwidth}
        \centering
        \includegraphics[width=\textwidth]{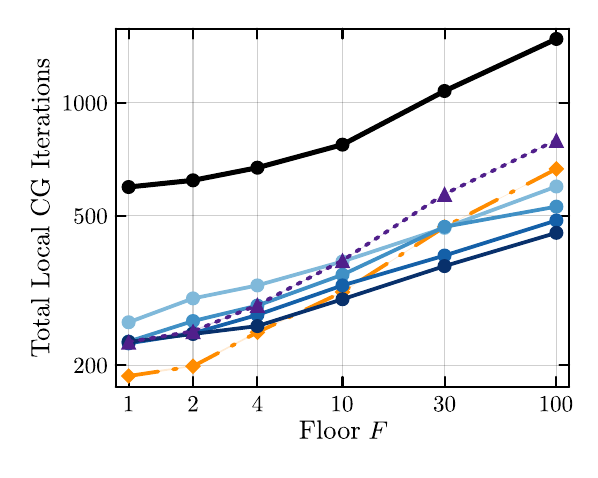}
        \caption{$N=20$}
        \label{fig:mals_iters_N20}
    \end{subfigure}
    \vspace{0.5em}
    \includegraphics[width=\textwidth]{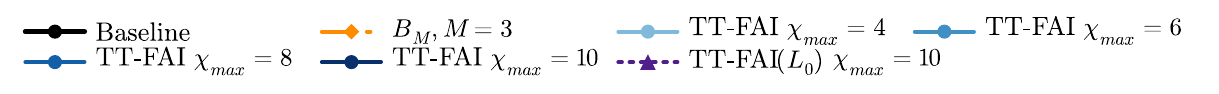}
    \caption{Total iterations of the local CG solver, summed over all sites and sweeps, within MALS to solve the
    anisotropic diffusion problem of Eq.~\eqref{eq:anisotropic diffusion} with $f \equiv 1$ up to a relative residual
    of $10^{-5}$, as a function of the floor $F$. Markers show the median over 10 repeats and the shaded bands the
    interquartile range. The $B_M$ curve corresponds to $M = 3$, since all the considered values of $M$ lie on the
    floor (see text).}
    \label{fig:mals_iters}
\end{figure}

Fig.~\ref{fig:mals_iters} shows the total number of local CG iterations. All the preconditioners reduce the
iterations of the baseline over the whole range of $F$, and their relative behavior follows the expectations
discussed above. For every $F$, however, the iterations of $B_M$ for all $M \in \{3,4,5,7,10\}$ differ by less than 7\% since it already reaches the floor with $M = 3$, at a bond dimension of 7, due to the coarseness of the grid.
For this reason, only $M = 3$ is shown in the figures. From the figure, we can see that as $F$ grows and the operator becomes less 
separable, the iterations of $B_M$ and the TT-FAI of $\mathcal{L}_0$ increase with the same trend, both subject to the same
floor, whereas the TT-FAI of $\mathcal{L}(\rho)$ falls below them from $F = 4$ ($\chi \geq 6$) for $N = 10$ and from $F = 10$ ($\chi \geq 10$) for $N = 20$.
Since the TT-FAI of $\mathcal{L}_0$ behaves like
$B_M$, the advantage of the TT-FAI of $\mathcal{L}(\rho)$ comes from capturing the coupling between dimensions rather than
from approximating $\mathcal{L}_0^{-1}$ better. At $N = 20$, however, the TT-FAI of $\mathcal{L}_0$ requires about 20\% more iterations than $B_M$ for every $F$. In this
case the build reaches its residual tolerance of $10^{-1}$ after a single sweep with a $G_L$ of bond dimension 6--7,
below the allowed $\chi_{\max} = 10$, so this TT-FAI is a coarser approximation of $\mathcal{L}_0^{-1}$ than at $N = 10$.
The same gap appears for $F = 1$, where both TT-FAI variants coincide, which confirms that it stems from the accuracy of the
TT-FAI construction and not from the coupling between dimensions.

\begin{figure}[htbp]
    \centering
    \begin{subfigure}[b]{0.49\textwidth}
        \centering
        \includegraphics[width=\textwidth]{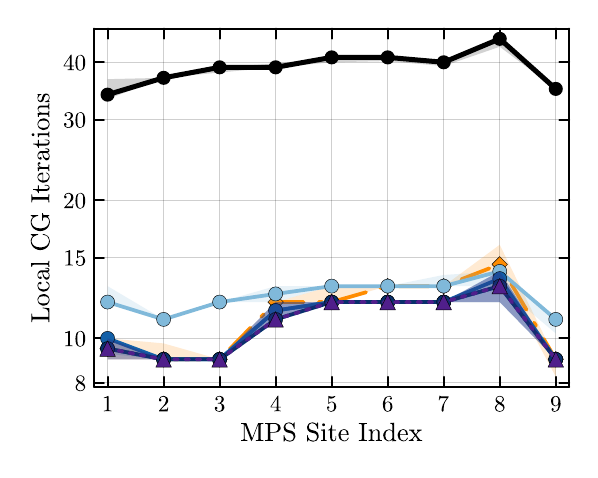}
        \caption{$N=10$, $F=1$}
        \label{fig:mals_site_iters_F1}
    \end{subfigure}\hfill
    \begin{subfigure}[b]{0.49\textwidth}
        \centering
        \includegraphics[width=\textwidth]{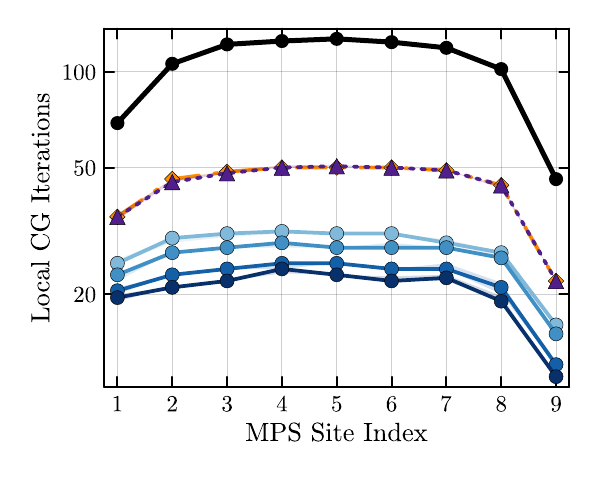}
        \caption{$N=10$, $F=100$}
        \label{fig:mals_site_iters_F100}
    \end{subfigure}
    \\[0.5em]
    \begin{subfigure}[b]{0.49\textwidth}
        \centering
        \includegraphics[width=\textwidth]{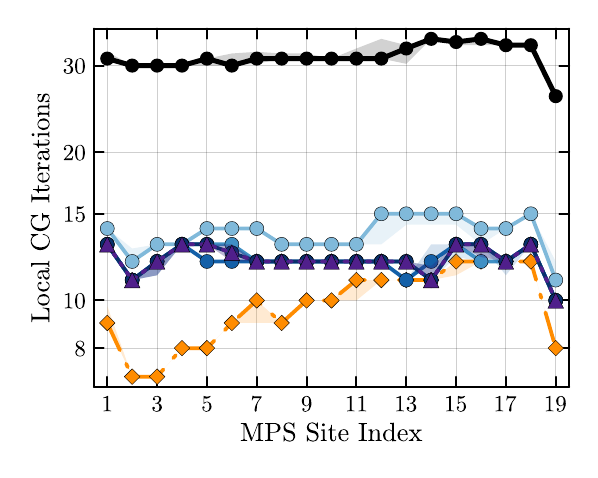}
        \caption{$N=20$, $F=1$}
        \label{fig:mals_site_iters_N20_F1}
    \end{subfigure}\hfill
    \begin{subfigure}[b]{0.49\textwidth}
        \centering
        \includegraphics[width=\textwidth]{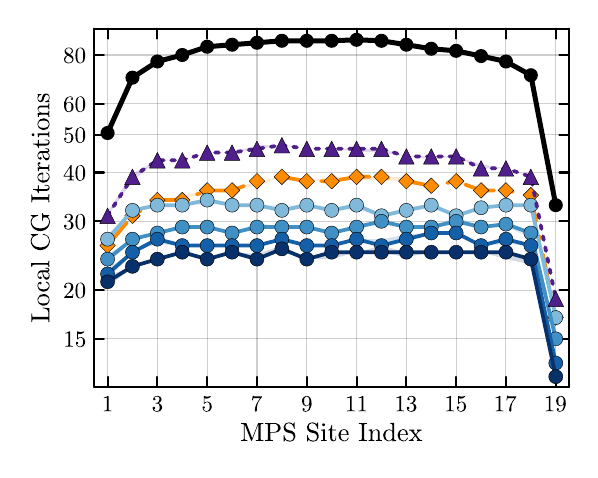}
        \caption{$N=20$, $F=100$}
        \label{fig:mals_site_iters_N20_F100}
    \end{subfigure}
    \vspace{0.5em}
    \includegraphics[width=\textwidth]{results/tt_anisotropy_stats/plots_paper_mugersh/pairs/legend_strip.pdf}
    \caption{Local CG iterations of MALS at each two-site block of the TT, summed over all sweeps, for the
    anisotropic diffusion problem with $N = 10$ (top) and $N = 20$ (bottom), $f \equiv 1$ and a relative tolerance of
    $10^{-5}$. Markers show the median over 10 repeats and the shaded bands the interquartile range. In all cases, the
    $B_M$ curve corresponds to $M = 3$.}
    \label{fig:mals_site_iters}
\end{figure}

In terms of site-by-site performance, Fig.~\ref{fig:mals_site_iters} shows how the local CG iterations are distributed 
along the TT for $F = 1$ and $F = 100$ for the two system sizes. For the separable case, the baseline requires 
around 30--45 iterations away from the boundary nodes, while
all the preconditioners bring this number down to around 10. For $F = 100$, at the center of the chain, the TT-FAI now clearly has an advantage,
with all the different bond dimension variants requiring around 20--30 iterations against the $\approx 50$ of $B_M$ and the
TT-FAI of $\mathcal{L}_0$ or the $\approx 125$ of the baseline at $N=10$.
At $N = 20$, the difference between $B_M$ and TT-FAI is less noticeable. Still, this confirms, as in Fig.~\ref{fig:mals_iters_N10}, 
that as the problem at hand becomes less separable, the quality of the $B_M$ preconditioners
starts to degrade and that the operator-agnostic approach of the TT-FAI is able to
better capture the correlations between the different terms of $\mathcal{L}$.

The reduction in iterations, however, does not translate into a reduction of the solve time, for any of the preconditioners considered. Excluding their construction, the median solve time stays within $0.94$--$1.40$ times that of the baseline for the TT-FAI of $\mathcal{L}(\rho)$ over all the values of $N$, $F$ and $\chi$, within $1.03$--$1.97$ for the TT-FAI of $\mathcal{L}_0$ and within $0.97$--$1.19$ for $B_3$, while larger values of $M$ are only slower.
This follows from the cost balance of Sec.~\ref{sec:Preconditioning ALS}. Since $\mathcal{L}$ has a bond dimension of
just three (Eq.~\eqref{eq:anisotropic mpo}), a preconditioned local iteration costs up to $1 + 2\chi/3 \approx 3.7$--$7.7$
unpreconditioned ones for the TT-FAI with $\chi = 4$--$10$, and $1 + 7/3 \approx 3.3$ for $B_3$, comparable to the iteration reductions of about $2$--$5$ observed across the
benchmark. This, alongside the increased cost given by the need to also keep track of and update the environments of the preconditioner, leads to the TT-FAIs of $\mathcal{L}(\rho)$ and $\mathcal{L}_0$, and $B_3$ being, on average, $11\%$, $30\%$ and $5\%$ slower than the baseline, respectively.

On top of this, the cost of actually building the preconditioner also needs to be taken into
account. For the $B_M$ approximate inverses, it is almost negligible, as their construction is direct from Eq.~\eqref{eq:sinc preconditioners}.
For the TT-FAI, however, its construction takes a median of about half a minute to 12 minutes depending on $N$, $F$ and $\chi$,
longer than the solve itself in every case. This time is mainly spent on solving Eq.~\eqref{eq:projected_fsai_linearsys_withI} on the vectorized space, in which the local basis size of $n^2=225$ for the
considered case makes the local iterations of the MALS sweeps very expensive. Taking this into account, therefore, and 
even though the results of the TT-FAI in terms of iteration reduction have been satisfactory in this setting, its use in 
practice is clearly limited.

\subsection{Two-dimensional Chorin projection}
\label{subsec:chorin_experiment}
In the previous experiment, the cost of building the TT-FAI preconditioner in a vectorized space outweighed its iteration savings. In this second experiment, we aim to tackle this issue by simulating a two-dimensional incompressible fluid in the QTT
format, replicating the numerical pipeline introduced in Ref.~\cite{peddinti2024e}. Fluid dynamics problems have been
studied in different TT settings \cite{Gourianov2022,peddinti2024e,Kiffner2023,hulst2026,kornev2024,danis2025}, but the use of
preconditioners has barely been discussed. This experiment is, nonetheless, well suited for approximate inversion
techniques such as the TT-FAI, as the QTT format, with a site dimension of just two, makes tackling
Eq.~\eqref{eq:projected_fsai_linearsys_withI} much less costly. On top of this, the Chorin projection scheme
\cite{chorin1967, chorin1968} employed solves a Poisson problem with the same operator at every time step, and,
therefore, the preconditioner can be constructed once and reused throughout the whole evolution of the fluid.

This pipeline proceeds in a two-step fashion. First, the fluid is evolved over a time step $dt$ without accounting for
the pressure, by integrating
\begin{equation}
    \label{eq:chorin predictor}
    \partial_t \mathbf{u} = \frac{1}{\mathrm{Re}} \nabla^2 \mathbf{u} - (\mathbf{u} \cdot \nabla) \mathbf{u}
\end{equation}
from $\mathbf{u}_n$ with an explicit scheme, which yields a provisional velocity field
$\mathbf{u}^*$. Then, a pressure field is obtained from $\mathbf{u}^*$ by imposing that $\mathbf{u}_{n+1}$
satisfies the incompressibility condition $\nabla \cdot \mathbf{u}_{n+1} = 0$, which results in the following Poisson
equation for the pressure
\begin{equation}
    \label{eq:pressure poisson}
    \nabla ^2 p_{n+1} = \frac{1}{dt}\nabla \cdot \mathbf{u}^*.
\end{equation}
Once a solution pressure field has been obtained, $\mathbf{u}^*$ is corrected to obtain the velocity field for
the next time step
\begin{equation}
    \mathbf{u}_{n+1} = \mathbf{u}^* - dt\, \nabla p_{n+1}.
\end{equation}
Note that, regardless of the order of the integrator used for Eq.~\eqref{eq:chorin predictor}, this splitting is only
first-order accurate in $dt$.

In order to tackle these equations, a finite difference scheme on a collocated grid is used to discretize the domain.
Within this setting, the $\nabla^2$ operator that enforces exactly the discrete incompressibility constraint in
Eq.~\eqref{eq:pressure poisson} is not the compact Laplacian operator, but the product of the discrete divergence and gradient \cite{ferziger2020}. This is a wide-stencil operator that is in general not symmetric for the inlet-outlet boundary conditions
used here, as the two operators carry the closures of the velocity and
pressure boundary conditions, respectively. Following common practice in collocated schemes, we instead replace this composite construction
by the compact Laplacian, built as
the sum of one-dimensional central second-derivative stencils of the same order as the rest of the operators
involved in the time stepping scheme (eighth order in our experiments), as in Ref.~\cite{peddinti2024e}. The resulting 
velocity field is therefore only approximately divergence-free, in the same way for all the solvers compared here, while 
the system becomes, after a symmetrization of the Neumann boundary rows, SPD with a low-rank QTT representation.

The simulated setup is a uniform flow past a cylinder in a channel. The domain is the box
$[0,2]\times[0,1]$, discretized with $2^{n_x}\times 2^{n_y}$ grid points.
The fields are encoded in a non-interleaved QTT, with the $n_x$ sites of the $x$ coordinate
followed by the $n_y$ sites of the $y$ coordinate, as in Ref.~\cite{peddinti2024e}. At the inlet ($x=0$), Dirichlet conditions $\mathbf{u}=(1,0)$ are
imposed on the velocity together with a homogeneous Neumann condition on the pressure, while at the outlet ($x=2$),
homogeneous Neumann conditions are imposed on the velocity and a Dirichlet condition with a constant reference value is
imposed on the pressure. The latter makes the pressure operator non-singular, so no regularization of the Poisson
problem is needed. The $y$ direction is periodic for all the fields. A cylinder of radius $0.1$ centered at $(0.5,0.5)$
is immersed in the fluid by multiplying the velocity field at every time step by the smoothed indicator function
\begin{equation}
    \text{mask}(x,y) = \frac{1}{2}\left[1 + \tanh\left(700\left((x-0.5)^2 + (y-0.5)^2 - 0.1^2\right)\right)\right],
\end{equation}
which is precomputed as an MPS with the Tensor Cross Interpolation algorithm \cite{savostyanov2011a} implementation
presented in Ref.~\cite{nunezfernandez2025}. Since the obstacle only enters through this mask on the
velocity, the pressure operator, and therefore its preconditioner, stays the same throughout the evolution. The fluid
is initialized with $\mathbf{u}_0 = (\text{mask}, 0)$ and evolved at $\mathrm{Re}=3000$, defined with respect to the inflow velocity
and the unit height of the channel, which results in an effective $\mathrm{Re}_D = 600$ based on the cylinder diameter $D=0.2$.

Regarding the comparison of the linear solvers applied to Eq.~\eqref{eq:pressure poisson}, the $B_M$ preconditioner of
Eq.~\eqref{eq:sinc preconditioners} cannot be applied here with the low bond dimension it has in the non-quantized setting. In Eq.~\eqref{eq:sinc preconditioners} each of the terms of the sum is of bond dimension one, given by 
each of the exponentials of the matrices $\mathcal{A}_\ell$ being promoted to the full space via identities on the remaining dimensions. In 
contrast, in the QTT setting studied here, the sum would involve exponentiating TT operators, which is usually
done via approximations which could make the operator indefinite, as shown in Fig.~\ref{fig:spd_nneg}, or working
with exponentially large matrices, which would defeat the purpose of the QTT formulation. We therefore compare the TT-FAI
against the baseline unpreconditioned solver only.

Table~\ref{tab:chorin_mals_params} lists the parameters used for the two MALS-based Chorin benchmark setups. The
predictor step of Eq.~\eqref{eq:chorin predictor} is integrated with a fourth-order Runge--Kutta scheme, and at each time
step the MALS sweeps are stopped once a relative residual of $10^{-3}$ is reached, using the residual-based truncation
of Ref.~\cite{oseledets2012} with no singular value cutoff, while the local CG solver is run up to a tighter tolerance to ensure the residual bound is reached. The TT-FAI is
built once per run and is then reused at every time step. Similarly to the anisotropy experiment, the Gershgorin bound has been used to determine $\mu$. In this case, the absolute values of the eighth-order stencil coefficients add up to
$\sigma = 205/72 + 2\left(8/5 + 1/5 + 8/315 + 1/560\right) = 2048/315 \approx 6.50$, so each dimension contributes
$\sigma/h_\ell^2$ to an interior row, resulting in
\begin{equation}
\label{eq:chorin_gershgorin}
    \mu = \|\Delta\|_\infty = \frac{2048}{315}\left(\frac{1}{h_x^2} + \frac{1}{h_y^2}\right),
\end{equation}
with $h_x = 2/(2^{n_x}-1)$ and $h_y = 1/2^{n_y}$ on the $[0,2]\times[0,1)$ box, which is periodic in $y$. At each step,
the pressure solve is initialized with the solution of the previous one. Every configuration has been repeated 10
times, each with a different random seed, and the figures below report the median over these
repeats. It is worth noting that the (10,10) system runs were not co-located, which adds node-to-node variation to the time measurements (17\% coefficient of variation for the baseline case) but not to the iteration counts.
\begin{table}[htbp]
    \centering
    \begin{tabular}{@{}cccccc@{}}
        \toprule
        $(n_x,n_y)$ & Grid & $\mu$ & $dt$ & Time steps & Max. bond dim. \\
        \midrule
        $(8,8)$ & $256^2$   & $5.3\times10^5$ & $10^{-3}$          & 500  & 150 \\
        $(10, 10)$ & $1024^2$ & $8.5\times10^6$ & $5\times10^{-4}$   & 1000 & 200 \\
        \bottomrule
    \end{tabular}
    \caption{MALS parameters for the Chorin benchmark, run at $\mathrm{Re}=3000$ with an eighth-order finite-difference
    stencil on a $2\times1$ box up to a final time of 0.5. The outer relative residual tolerance is $10^{-3}$ with 100 sweeps at most, with a residual-based truncation tolerance of $5 \times 10^{-4}$ and no singular value cutoff, so the truncation is controlled by the residual alone. The local CG solver runs up to a relative tolerance of $10^{-4}$ with at most 2000 iterations per local problem. The TT-FAI is built in all cases with at most 30 MALS sweeps at a local tolerance of $10^{-8}$, 3000 maximum CG iterations and a global residual tolerance of $10^{-1}$. }
    \label{tab:chorin_mals_params}
\end{table}

\begin{figure}[htbp]
    \centering
    \begin{subfigure}[b]{0.49\textwidth}
        \centering
        \includegraphics[width=\textwidth]{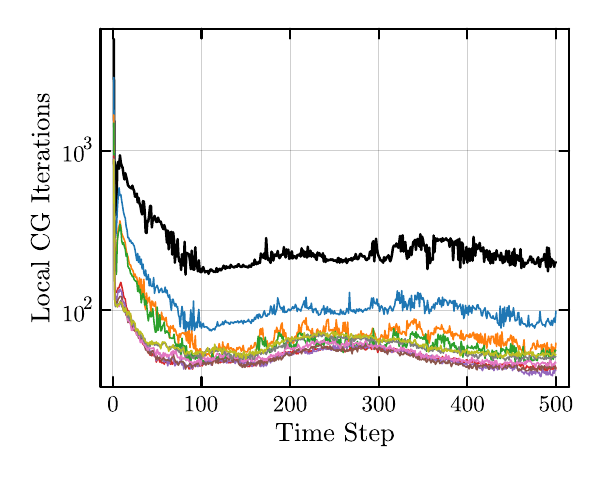}
        \caption{$(n_x, n_y) = (8,8)$, $256^2 \text{ grid}$}
 
    \end{subfigure}\hfill
    \begin{subfigure}[b]{0.49\textwidth}
        \centering
        \includegraphics[width=\textwidth]{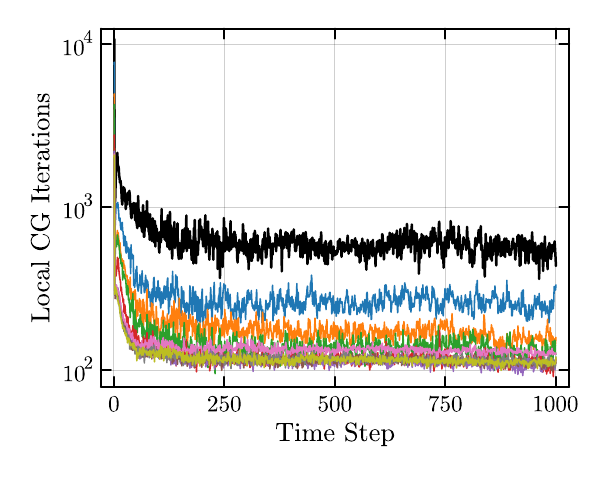}
        \caption{$(n_x, n_y) = (10,10)$, $1024^2 \text{ grid}$}
    \end{subfigure}
    \vspace{0.5em}
    \includegraphics[width=\textwidth]{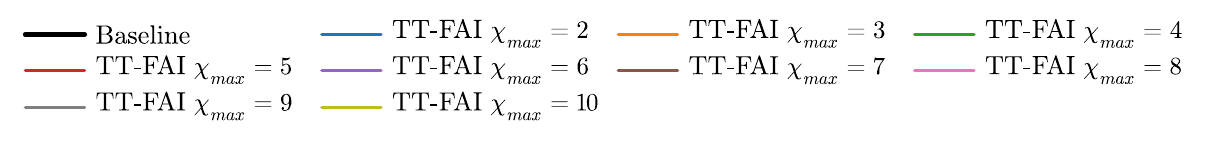}
    \caption{Total iterations of the local CG solver within MALS, summed over all sites and sweeps, to solve the
    Poisson problem of Eq.~\eqref{eq:pressure poisson} at each time step, over time evolutions of 500 steps for
    $(n_x, n_y) = (8,8)$ and 1000 steps for $(n_x, n_y) = (10,10)$. Curves show the median over 10 repeats.}
    \label{fig:chorin_mals_iterations}
\end{figure}

\begin{figure}[htbp]
    \centering
    \begin{subfigure}[b]{0.49\textwidth}
        \centering
        \includegraphics[width=\textwidth]{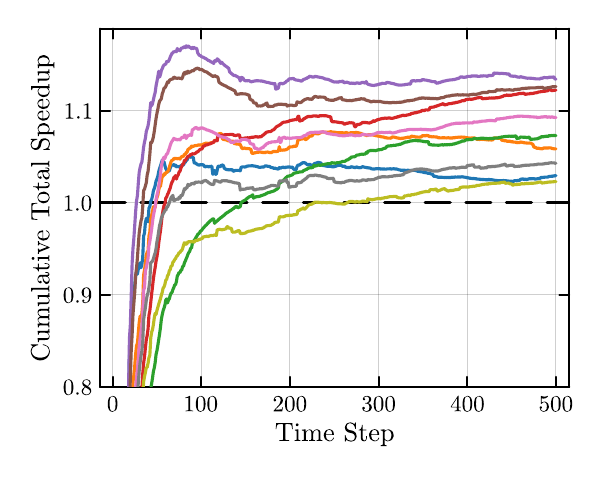}
        \caption{$(n_x, n_y) = (8,8)$, $256^2 \text{ grid}$}
    \end{subfigure}\hfill
    \begin{subfigure}[b]{0.49\textwidth}
        \centering
        \includegraphics[width=\textwidth]{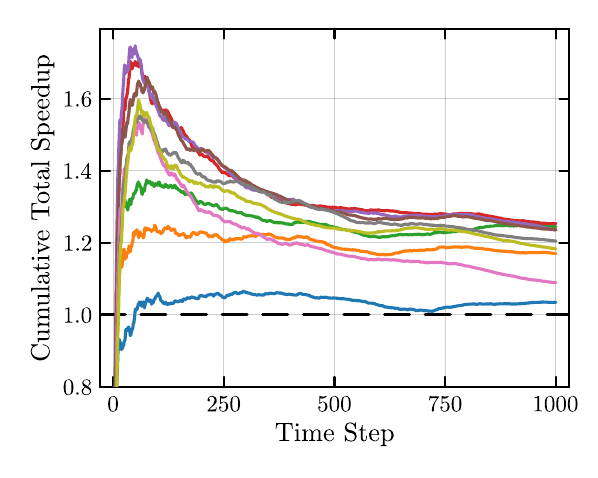}
        \caption{$(n_x, n_y) = (10,10)$, $1024^2 \text{ grid}$}
    \end{subfigure}
    \vspace{0.5em}
    \includegraphics[width=\textwidth]{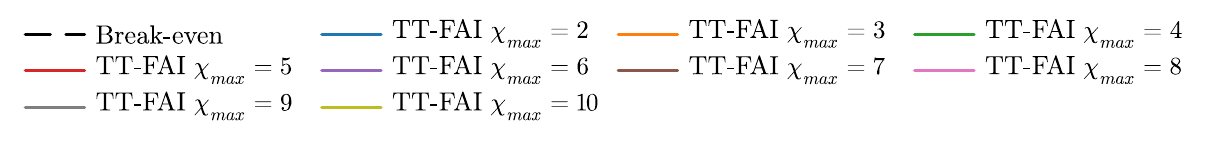}
    \caption{Cumulative total speedup, measured as the ratio between the times of the unpreconditioned MALS and of each TT-FAI variant to solve the Poisson problem of Eq.~\eqref{eq:pressure poisson} up to each time
    step, over time evolutions of 500 steps for $(n_x, n_y) = (8,8)$ and 1000 steps for $(n_x, n_y) = (10,10)$, taking into account also 
    preconditioner build time. The lines show the median of the ratio between baseline and TT-FAI runs over 10 repeats, each with a different random seed. In both panels, the
    vertical axis is cut at $0.8$ for readability, as the curves exceed $1$ after 7--13 steps at (10,10) (50 for the $\chi=2$ case) and after 27--265 steps at (8,8).}
    \label{fig:chorin_mals_cumulative_speedup}
\end{figure}

Fig.~\ref{fig:chorin_mals_iterations} shows the total number of local CG iterations for each time step of the evolution
of both system sizes. As in the anisotropy experiment considered in the previous subsection, we observe how the TT-FAI 
preconditioners successfully reduce the number of iterations, with larger iteration count differences for higher bond
dimensions up to $\chi = 6$. Also, this iteration reduction does provide an actual speedup, as shown in Fig.~\ref{fig:chorin_mals_cumulative_speedup},
in which all of the TT-FAI setups are able to provide a runtime improvement, especially for the more ill-conditioned system with
$10$ sites per dimension, in which a TT-FAI with a bond dimension between $4$ and $7$ is able to achieve a final net speedup of about 20--25\% at the end of the 1000-step run. For the smaller system size, since the operator
is not so ill-conditioned due to the coarse grid, the overall number of iterations performed by both the baseline and the 
TT-FAI preconditioned runs gets reduced. This makes the overhead given by the tracking and updating of the environments
that make up the projectors $U$ depicted in Fig.~\ref{fig:projectedGlTGl} larger in relative terms,
diluting the advantage given by the reduction in iterations. The same effect explains the decay in the speedup curves of the (10,10) system in Fig.~\ref{fig:chorin_mals_cumulative_speedup}b, as warm starts lower the iterations per step.

We attribute this net speedup to two main features of this 
linear system. The first one is that the boundary corrections made to the MPO over the baseline Laplacian construction 
from Ref.~\cite{Kazeev2012} result in a TT operator where most bonds have dimension 4, compared to the 3 in the 
previous experiment, which results in the relative overhead of the TT-FAI $2\chi_{G_L}/\chi_{A}$ being slightly lower. And, secondly, since the QTT formulation of the problem allows exploring 
an exponentially fine grid, the condition number of the corresponding operator grows exponentially with the number of sites, as opposed to the $N$-dimensional case studied in the previous subsection, hindering
the convergence of the CG method, which gives more room for the preconditioner to reduce the number of iterations. 
On top of this, since the local dimension of the QTT is two, the solution
of Eq.~\eqref{eq:projected_fsai_linearsys_withI} is much cheaper. This leads to the iterative scheme being able to 
amortize the build time of the preconditioner in just a few time steps, 
as can be seen from Fig.~\ref{fig:chorin_mals_cumulative_speedup}.

\section{Conclusions}
We have presented a numerical pipeline to obtain a factorized approximate inverse in the TT format (TT-FAI). This 
construction, unlike previous direct approximate inversion techniques, does not suffer from indefiniteness caused by 
truncation errors, which results in a suitable preconditioner for TT linear solvers based on the Conjugate Gradient
method. By means of an efficient projection step and a subsequent linear solve in a vectorized space, the factors of the 
preconditioner can be obtained variationally without the need
to compute the Cholesky decomposition of the linear system's operator, with a bond dimension that controls the cost of the approximate inverse at 
the expense of preconditioner quality. Within MALS, projecting the two factors separately 
keeps every local preconditioner positive semidefinite.

We have benchmarked the TT-FAI preconditioners in two different settings, both using such MALS solvers. First, we considered an $N$-dimensional rotated 
anisotropy linear system, which can be tuned to reduce to a Poisson equation. As the experiment is tuned away from this separable regime, the TT-FAI approach is able to capture the correlations between the different dimensions and reduce the iterations by up to 50\% for 10 dimensions and up to 30\% for 20 dimensions with respect to the preconditioner of Ref.~\cite{khoromskij2009b}.
With the rotated anisotropy having a bond dimension of just 3, however, each preconditioned iteration costs up to 3.7--7.7 unpreconditioned ones, resulting in an actual net slowdown. On top of this, the build takes several times longer than the actual solve, because of the local dimension $n^2 = 225$
of the vectorized problem, which rules out the use of TT-FAI for practical purposes in such a setting.

As a second benchmark we have studied the time evolution of a two-dimensional incompressible fluid by means of a Chorin
projection scheme and a QTT formulation. This setting has three main advantages. The site dimension is two, which
keeps the build cheap, the operator of the Poisson linear system has a larger bond dimension than in the anisotropy case (4--6 against 3), and a single TT-FAI is reused at every time step. On a $1024^2$ grid, TT-FAIs
with $\chi = 4$--$7$ reduce the iterations by a factor of 4.0--5.6 and give net speedups of about 20--25\% over 1000
time steps, build included. On the coarser $256^2$ grid the solves need fewer iterations, so the fixed work of each
MALS step dilutes the net speedup to at most 13\%.

The TT-FAI is therefore a practical addition to TT solvers when the local dimension is small and the preconditioner can
be reused, as in time stepping schemes. Still, two main questions remain open for future work. 
The first is how well low-rank TT-FAIs approximate the inverse of less structured operators, in more complex
geometries, for example, or within unstructured grids in the QTT format \cite{markeeva2021,kornev2024}. The second,
also related to such unstructured systems, is that the diagonal rescaling of the inverse factors requires the computed 
diagonal of the factor to stay positive. Lemma~\ref{lem:diag_positivity} shows that this can
only fail if the error of a whole row concentrates on it, which the nearly constant diagonal of uniform grids makes
unlikely, but the argument weakens with less structured operators. Benchmarking the TT-FAI in such settings would therefore
test the bounds given and also the reliability and cost of the energy-based rescaling of Sec.~\ref{subsec:computing_D}
in practical scenarios.

\section*{CRediT authorship contribution statement}
\textbf{Pau Torrente Badia:} Conceptualization, Methodology, Software, Formal analysis, Investigation, Writing -- Original Draft; \textbf{Daniel Pérez García:} Software, Validation, Writing -- Review \& Editing; \textbf{Artur Garcia-Saez:} Writing -- Review \& Editing, Supervision, Funding Acquisition; \textbf{Stefano Carignano:} Writing -- Review \& Editing, Supervision, Funding Acquisition.

\section*{Code and data availability}
The methods of this work are implemented in the open-source Julia package DifferenTT.jl, available under the MIT
license at \url{https://github.com/bsc-quantic/DifferenTT.jl} (version 0.1.0). The folder \texttt{experiments/paper}
of the repository contains the scripts that produce every figure and table with numerical results in this paper,
together with instructions to rerun them. The data shown in these figures and tables were produced at the tag
\texttt{paper-ttfai-v0.1-data} of the repository, whose library code is that of version 0.1.0, and are archived,
with the Julia environment used, on Zenodo under the CC BY 4.0 license: \url{https://doi.org/10.5281/zenodo.23663326}.

\section*{Declaration of competing interest}
The authors declare that they have no known competing financial interests or personal relationships that could have appeared to influence the work reported in this paper.

\section*{Acknowledgments}
We thank Aleix Bou-Comas for fruitful discussions on TT-based numerical methods and Guillaume Houzeaux, Ramiro de Olazábal and Mariano Vázquez for valuable insights on matrix-based linear system solvers.
We acknowledge funding from the project PID2023-151384NB-I00 
funded by MICIU/AEI/10.13039/501100011033 and by ERDF, EU.
S.C. acknowledges his AI4S fellowship
within the “Generación D” initiative by Red.es, Ministerio para la Transformación Digital y de la Función Pública, for talent attraction (C005/24-ED CV1),
funded by NextGenerationEU through PRTR. 

\section*{Declaration of generative AI and AI-assisted technologies in the manuscript preparation process}
During the preparation of this work the authors used Claude (Anthropic) and Gemini (Google) LLMs to assist in the formalization of Lemma~\ref{lem:diag_positivity} and the rotated anisotropy benchmark of Sec.~\ref{subsec:anisotropy_experiment}, as well as in formatting and revising the code used to perform the numerical experiments. After using these tools, the authors reviewed and edited the content as needed and take full responsibility for the content of the article.

\bibliographystyle{elsarticle-num}
\bibliography{references/applications,references/classicalsolvers,references/preconditioners,references/quantum,references/ttops,references/ttsolvers}

\newpage

\appendix

\section{Positivity of the diagonal of \texorpdfstring{$\tilde{G}_L$}{G L}}
\label{app:diag_positivity}

In this appendix we prove Lemma~\ref{lem:diag_positivity} and Corollary~\ref{cor:diag_fraction}, and discuss why
approximate solutions of Eq.~\eqref{eq:projected_fsai_linearsys_withI} do not produce non-positive diagonal entries in
practice.

The argument extends the classical one of Kolotilina and Yeremin \cite{kolotilina1993} to inexact solutions. There,
the scaling is defined through the energy of each row, $(\tilde{G}_LA\tilde{G}_L^T)_{ii}$, and replaced by
$(\tilde{G}_L)_{ii}$ through an identity that holds only when every row is computed exactly. The positivity of the
diagonal then follows from that of the diagonal of the inverse of an SPD matrix \cite[Lemma 2.1]{kolotilina1993}. Part
(a) of Lemma~\ref{lem:diag_positivity} is this classical result, written for the dense lower triangular pattern. Parts
(b) and (c) are what changes when the rows are only approximate, as is the case for any TT solver. They quantify how far
$(\tilde{G}_L)_{ii}$ can drift from its exact value, and how inaccurate a row has to be before its diagonal entry
changes sign.

\subsection{The build problem splits into one problem per row}
\label{app:row-wise splitting}
Let $A \in \mathbb{R}^{m\times m}$ be SPD, and let $\mathcal{T} := \mathcal{P}_L(I\otimes A)\mathcal{P}_L +
\mu(\mathcal{I} - \mathcal{P}_L)$ with $\mu > 0$. Consistently with Eq.~\eqref{eq:fsai_linearsystem}, $\text{vec}$ stacks
the rows of its argument, $\mathbf{x} = \text{vec}(\tilde{G}_L) = (\mathbf{x}_1;\dots;\mathbf{x}_m)$, where
$\mathbf{x}_i \in \mathbb{R}^m$ is the $i$-th row of $\tilde{G}_L$ written as a column vector, and $x_{ij} =
(\tilde{G}_L)_{ij}$ denotes its entries. In this ordering, $I\otimes A$ applies $A$ to every row
separately, and $\mathcal{P}_L$ keeps the entries $i \ge j$ of row $i$ and sets the others to zero. Neither operator mixes
different rows, so $\mathcal{T}$ is block diagonal, $\mathcal{T} = \text{blockdiag}(\mathcal{T}_1,\dots,\mathcal{T}_m)$.
Splitting each row into its lower entries $i \ge j$ and its upper entries $i < j$, the block of row $i$ reads
\begin{equation}
    \label{eq:app row block}
    \mathcal{T}_i = \begin{pmatrix} A_i & 0 \\ 0 & \mu\,\mathcal{I}_{m-i} \end{pmatrix},
    \qquad A_i = A[1{:}i,\,1{:}i],
\end{equation}
where $A_i$ is the leading $i\times i$ submatrix of $A$. $A_i$ is SPD because $A$ is, and therefore so is
$\mathcal{T}_i$. The right-hand side $\text{vec}(I)$ restricted to row $i$ is the canonical vector
$\mathbf{e}_i$, which lies in the lower part of the row. Eq.~\eqref{eq:projected_fsai_linearsys_withI} thus
decouples into $m$ independent systems,
\begin{equation}
    \label{eq:app row systems}
    \mathcal{T}_i\,\mathbf{x}_i^\star = \mathbf{e}_i, \qquad i = 1,\dots,m,
\end{equation}
whose solutions have vanishing upper entries and lower entries given by $A_i^{-1}\mathbf{e}_i$. This is exactly the
classical FSAI row condition of Eq.~\eqref{eq:fsai_linearsystem}, independently of $\mu$.

We equip the space of row $i$ with the energy inner product and norm
\begin{equation}
    \langle \mathbf{v}, \mathbf{w}\rangle_{\mathcal{T}_i} = \mathbf{v}^T\mathcal{T}_i\mathbf{w}, \qquad
    \|\mathbf{v}\|_{\mathcal{T}_i}^2 = \langle \mathbf{v}, \mathbf{v}\rangle_{\mathcal{T}_i}.
\end{equation}
Since $\mathcal{T}$ is SPD, Eq.~\eqref{eq:projected_fsai_linearsys_withI} is the stationarity condition of the
quadratic functional $J(\mathbf{x}) = \frac{1}{2}\mathbf{x}^T\mathcal{T}\mathbf{x} - \text{vec}(I)^T\mathbf{x}$,
which splits into a sum of row energies,
\begin{equation}
    \label{eq:app row energy}
    J(\mathbf{x}) = \sum_{i=1}^m J_i(\mathbf{x}_i), \qquad
    J_i(\mathbf{x}_i) = \frac{1}{2}\|\mathbf{x}_i\|_{\mathcal{T}_i}^2 - x_{ii}.
\end{equation}
Completing the square with Eq.~\eqref{eq:app row systems} gives
\begin{equation}
    \label{eq:app energy error}
    J_i(\mathbf{x}_i) - J_i(\mathbf{x}_i^\star) = \frac{1}{2}\|\mathbf{x}_i - \mathbf{x}_i^\star\|_{\mathcal{T}_i}^2,
\end{equation}
so the energy measures the distance to the exact row. We stress that Eqs.~\eqref{eq:app row energy} and~\eqref{eq:app energy error} hold for any $\mathbf{x}$. The TT format restricts which $\mathbf{x}$ a solver can reach, but not how $J$ evaluates it, so everything below applies to the output of any TT solver.

\subsection{The diagonal entry as an energy inner product}
Using Eq.~\eqref{eq:app row systems} to replace $\mathbf{e}_i$ by $\mathcal{T}_i\mathbf{x}_i^\star$, the diagonal entry
of any row $\mathbf{x}_i$ can be written as
\begin{equation}
    \label{eq:app diag inner prod}
    x_{ii} = \mathbf{e}_i^T\mathbf{x}_i = (\mathcal{T}_i\mathbf{x}_i^\star)^T\mathbf{x}_i
    = \langle \mathbf{x}_i^\star, \mathbf{x}_i\rangle_{\mathcal{T}_i}.
\end{equation}
In words, the diagonal entry of a row is its overlap with the exact row, measured in the energy inner product.
Its sign is therefore a statement about angles: $x_{ii} > 0$ if and only if $\mathbf{x}_i$ forms an angle smaller
than $90^\circ$ with $\mathbf{x}_i^\star$ in the geometry defined by $\mathcal{T}_i$. Two consequences follow directly.
Setting $\mathbf{x}_i = \mathbf{x}_i^\star$,
\begin{equation}
    \label{eq:app exact diag norm}
    g_i^\star = x_{ii}^\star = \|\mathbf{x}_i^\star\|_{\mathcal{T}_i}^2,
\end{equation}
so the exact diagonal entry is a squared norm. And, expanding the error of an arbitrary row,
\begin{equation}
    \label{eq:app error expansion}
    \|\mathbf{x}_i - \mathbf{x}_i^\star\|_{\mathcal{T}_i}^2 = \|\mathbf{x}_i\|_{\mathcal{T}_i}^2 - 2x_{ii} + g_i^\star.
\end{equation}

Eq.~\eqref{eq:app exact diag norm} is the identity behind the shortcut of Eq.~\eqref{eq:fsai_rescaling}. For a lower
triangular row, $\|\mathbf{x}_i\|_{\mathcal{T}_i}^2 = (\tilde{G}_LA\tilde{G}_L^T)_{ii}$ is the energy of the row, used by
Kolotilina and Yeremin to define the scaling, while $x_{ii} = \langle \mathbf{x}_i^\star,
\mathbf{x}_i\rangle_{\mathcal{T}_i}$ is what Eq.~\eqref{eq:rescaling D} actually uses. The two coincide at the exact
solution. Away from it, the energy remains positive for every non-zero row, which is why the iterative FSAI variants of
\cite{janna2015} scale with it, whereas the sign of $x_{ii}$ depends on how close $\mathbf{x}_i$ is to
$\mathbf{x}_i^\star$. Lemma~\ref{lem:diag_positivity} makes this dependence precise.

\subsection{Proof of Lemma~\ref{lem:diag_positivity}}

\paragraph{Part (a)} By Eq.~\eqref{eq:app row systems}, $g_i^\star = \mathbf{e}_i^T A_i^{-1}\mathbf{e}_i =
(A_i^{-1})_{ii}$, which is positive since $A_i^{-1}$ is SPD. This is the argument of
\cite[Lemma 2.1]{kolotilina1993}, where $A_i$ is replaced by the principal submatrix of $A$ selected by the sparsity
pattern of row $i$. To relate it to the Cholesky factor, note that the
leading $i\times i$ block of $L_A$, denoted $L_i$, is the Cholesky factor of $A_i$, so that $A_i^{-1} = L_i^{-T}L_i^{-1}$
and
\begin{equation}
    (A_i^{-1})_{ii} = \sum_{k=1}^{i} \left[(L_i^{-1})_{ki}\right]^2 = \left[(L_i^{-1})_{ii}\right]^2 =
    \frac{1}{(L_A)_{ii}^2},
\end{equation}
where only $k = i$ survives because $L_i^{-1}$ is lower triangular. Regarding the third equality, $(L_i^{-1})_{ii} = 1/(L_A)_{ii}$ because,
$L_i$ and $L_i^{-1}$ being both lower triangular, the $i$-th diagonal entry of $L_iL_i^{-1} = I$ reduces to
$(L_i)_{ii}(L_i^{-1})_{ii} = 1$, with $(L_i)_{ii} = (L_A)_{ii}$. For the lower bound, the $i$-th diagonal entry of
$A = L_AL_A^T$ is $A_{ii} = \sum_{k\le i}(L_A)_{ik}^2 \ge (L_A)_{ii}^2$, hence $g_i^\star \ge 1/A_{ii}$. For the upper
bound, $g_i^\star = \mathbf{e}_i^TA_i^{-1}\mathbf{e}_i \le \lambda_{\max}(A_i^{-1}) = 1/\lambda_{\min}(A_i)$, and
$\lambda_{\min}(A_i) \ge \lambda_{\min}(A)$ by Cauchy's interlacing theorem.

As a consistency check, Eq.~\eqref{eq:rescaling D} applied to the exact solution gives $D_{ii} = (g_i^\star)^{-1/2} =
(L_A)_{ii}$, so that the diagonal of $G_L = D\tilde{G}_L^\star$ is $(L_A)_{ii}/(L_A)_{ii}^2 = 1/(L_A)_{ii} =
(L_A^{-1})_{ii}$, as expected for $G_L = L_A^{-1}$.

\paragraph{Part (b)} Subtracting Eq.~\eqref{eq:app exact diag norm} from Eq.~\eqref{eq:app diag inner prod} and applying
the Cauchy--Schwarz inequality in the energy inner product,
\begin{equation}
    |x_{ii} - g_i^\star| = |\langle \mathbf{x}_i^\star, \mathbf{x}_i - \mathbf{x}_i^\star\rangle_{\mathcal{T}_i}|
    \le \|\mathbf{x}_i^\star\|_{\mathcal{T}_i}\,\|\mathbf{x}_i - \mathbf{x}_i^\star\|_{\mathcal{T}_i}.
\end{equation}
Dividing by $g_i^\star = \|\mathbf{x}_i^\star\|_{\mathcal{T}_i}^2$ gives $|x_{ii}/g_i^\star - 1| \le \nu_i$.

\paragraph{Part (c)} If $x_{ii} \le 0$, then $x_{ii}/g_i^\star - 1 \le -1$, and part (b) gives $\nu_i \ge 1$. The
same conclusion follows from Eq.~\eqref{eq:app error expansion}, which also shows when the bound is attained: for
$x_{ii} \le 0$ all three terms on its right-hand side are non-negative, so
$\|\mathbf{x}_i - \mathbf{x}_i^\star\|_{\mathcal{T}_i}^2 \ge g_i^\star = \|\mathbf{x}_i^\star\|_{\mathcal{T}_i}^2$,
with equality only for $\mathbf{x}_i = 0$. Among all rows with a non-positive diagonal entry, the one closest to the
exact row is the zero row: the cheapest way of changing the sign of the diagonal is to discard the whole row. \hfill$\square$

\subsection{Proof of Corollary~\ref{cor:diag_fraction}}
Let $\mathcal{N} = \{i : x_{ii} \le 0\}$. By part (c), every $i \in \mathcal{N}$ satisfies
$\|\mathbf{x}_i - \mathbf{x}_i^\star\|_{\mathcal{T}_i}^2 \ge g_i^\star$. Since $\mathcal{T}$ is block diagonal,
\begin{equation}
    \sum_{i\in\mathcal{N}} g_i^\star \le \sum_{i\in\mathcal{N}} \|\mathbf{x}_i - \mathbf{x}_i^\star\|_{\mathcal{T}_i}^2
    \le \sum_{i=1}^m \|\mathbf{x}_i - \mathbf{x}_i^\star\|_{\mathcal{T}_i}^2 = \|\mathbf{x} - \mathbf{x}^\star\|_{\mathcal{T}}^2,
\end{equation}
while, by Eq.~\eqref{eq:app exact diag norm}, $\|\mathbf{x}^\star\|_{\mathcal{T}}^2 = \sum_i g_i^\star =
\text{Tr}(\tilde{G}_L^\star)$. Dividing the two expressions gives Eq.~\eqref{eq:diag fraction}. \hfill$\square$

\subsection{Why sign changes are not observed in practice}
\label{app:sign_changes}
Lemma~\ref{lem:diag_positivity} and Corollary~\ref{cor:diag_fraction} turn the question into one about how the error of
the computed $\tilde{G}_L$ is distributed over the rows. Three observations explain why no row reaches $\nu_i \ge 1$.

\paragraph{The energy budget} By Corollary~\ref{cor:diag_fraction}, a global relative energy error of $\varepsilon_E =
10^{-2}$ allows at most a weighted fraction $10^{-4}$ of the rows to change sign, and each of them must be completely
wrong. A sign change is therefore never a small perturbation, as it requires concentrating an error of the size of a
whole row on that row.

\paragraph{The exact diagonal is nearly constant} By part (c), row $i$ can only change sign if its error energy
reaches $g_i^\star$, so $g_i^\star$ is the price of a sign change in that row. If these prices differed by orders of
magnitude, the rows with the smallest $g_i^\star$ would be the first to flip, and the weighted fraction of
Eq.~\eqref{eq:diag fraction} could hide many of them. For the operators considered here, the prices are all alike.
By part (a), $g_i^\star = 1/(L_A)_{ii}^2$, where $(L_A)_{ii}^2$ is what remains of $A_{ii}$ after eliminating the
unknowns that precede $i$ in the ordering. On a uniform grid with constant coefficients, $A_{ii}$ is the same for every row and, for a local stencil, these remainders settle to a common value outside thin layers next to the boundaries.
No row is therefore markedly cheaper to flip than the others, and the weighted fraction of
Eq.~\eqref{eq:diag fraction} becomes, up to a factor of order one, a plain count of rows. This no longer holds on
strongly non-uniform grids or with strongly varying coefficients, as discussed in \ref{app:diag_limits}.

\paragraph{A flat diagonal is easy to represent in the TT format} A nearly constant $\text{diag}(\tilde{G}_L^\star)$,
with thin boundary layers, has TT rank close to one, and is among the first features that any low-rank solver captures.
To change the sign of a single entry, the error of the computed $\tilde{G}_L$ would need a spike of the size of the
whole row localized on that row, while remaining small elsewhere. Such a localized feature has a high TT rank, which
neither a rank truncation nor an incomplete convergence of the solver has a reason to produce: their errors are spread
over many rows, so that each row sees $\nu_i \approx \varepsilon_E \ll 1$. Unlike the previous two observations, this
one is a heuristic and not a theorem.

\subsection{What the lemma does not guarantee}
\label{app:diag_limits}
\paragraph{The stopping criterion does not control $\varepsilon_E$} TT solvers monitor the relative residual
$\|\mathbf{r}\|/\|\mathbf{b}\|$, with $\mathbf{r} = \mathbf{b} - \mathcal{T}\mathbf{x}$ and $\mathbf{b} =
\text{vec}(I)$, rather than the energy error. Since $\|\mathbf{x} - \mathbf{x}^\star\|_{\mathcal{T}}^2 =
\mathbf{r}^T\mathcal{T}^{-1}\mathbf{r} \le \|\mathbf{r}\|^2/\lambda_{\min}(\mathcal{T})$ and
$\|\mathbf{x}^\star\|_{\mathcal{T}}^2 = \mathbf{b}^T\mathcal{T}^{-1}\mathbf{b} \ge
\|\mathbf{b}\|^2/\lambda_{\max}(\mathcal{T})$, the two are related, in the worst case, by
\begin{equation}
    \varepsilon_E \le \sqrt{\kappa(\mathcal{T})}\,\frac{\|\mathbf{r}\|}{\|\mathbf{b}\|},
\end{equation}
which carries no information for a build tolerance of $10^{-1}$ and an ill-conditioned $A$.

\paragraph{A bound on a fraction is not a bound on existence} Eq.~\eqref{eq:diag fraction} limits the share of rows
with a non-positive diagonal, but it does not exclude them. For the problem sizes considered here, with up to
$15^{20}$ rows, even a tiny fraction amounts to many entries if the error concentrates. The heuristic above explains why
it does not, but it is not a proof.

\paragraph{Non-uniform grids and varying coefficients} The flat diagonal discussed above relies on all
rows of $A$ having a similar diagonal entry. On a non-uniform grid this is not the case: a finite-difference diagonal
grows as $h^{-2}$ with the local mesh size $h$, so refining the mesh by a factor of ten in some region makes $A_{ii}$
there a hundred times larger, and $g_i^\star$ a hundred times smaller. The same happens where the coefficients of the
PDE are large. These rows are the cheapest to flip, and therefore the first to change sign if the error concentrates
on them. In addition, the diagonal of $\tilde{G}_L^\star$ is no longer nearly constant, so its TT rank grows and it is
harder for the solver to capture accurately.

\paragraph{A practical diagnostic} Since $g_i^\star A_{ii} \ge 1$ by part (a), sampled values of
$(\tilde{G}_L)_{ii}A_{ii}$ noticeably below one indicate that the computed diagonal underestimates the exact one in
those rows, leading to too large entries of $D$, and non-positive values indicate a sign change. Each sample requires
evaluating a single entry of $\tilde{G}_L$ and of $A$ in the TT format.

\section{Derivation of the anisotropy floor}
\label{app:anisotropy_floor}
In this appendix we derive, step by step, the bound of Eq.~\eqref{eq:anisotropy floor}, which is the quantity used in
Sec.~\ref{subsec:anisotropy_experiment} to control the difficulty of the anisotropic diffusion benchmark for separable
preconditioners. As in~\ref{app:diag_positivity}, instead of comparing the operators $\mathcal{L}(\rho)$ and $\mathcal{L}_0$ directly, the argument relies on comparing the energies $\mathbf{u}^T\mathcal{L}(\rho)\,\mathbf{u}$ and $\mathbf{u}^T\mathcal{L}_0\,\mathbf{u}$ that they assign to an arbitrary grid vector $\mathbf{u}$, which reduces every
step of the derivation to an inequality between real numbers.

\subsection{Notation and preliminaries}
We consider $n$ grid points per dimension with spacing $h$ and homogeneous Dirichlet boundary
conditions at the ghost nodes, so that, for a one-dimensional vector $\mathbf{u} \in \mathbb{R}^n$ extended with $u_0 = u_{n+1} = 0$, the
second and first difference matrices of Eq.~\eqref{eq:anisotropic mpo} act as
\begin{equation}
    \label{eq:app stencils}
    (\mathcal{A}\mathbf{u})_i = \frac{-u_{i-1} + 2u_i - u_{i+1}}{h^2}, \qquad
    (\mathcal{C}\mathbf{u})_i = \frac{u_{i+1} - u_{i-1}}{2h}, \qquad i = 1, \dots, n.
\end{equation}
In $N$ dimensions, $\mathcal{A}_\ell$ and $\mathcal{C}_\ell$ denote these matrices acting along the $\ell$-th dimension,
i.e. the Kronecker product of $\mathcal{A}$ or $\mathcal{C}$ in the $\ell$-th position with identities in all the
others, and $\mathcal{L}_0 = \sum_\ell \mathcal{A}_\ell = -\Delta_h$.

For two symmetric matrices $X$ and $Y$, we write $X \preceq Y$ whenever
$\mathbf{u}^TX\mathbf{u} \le \mathbf{u}^TY\mathbf{u}$ for every vector $\mathbf{u}$, or equivalently, whenever $Y - X$ is
positive semidefinite. The bound we are after concerns the spectrum of $\mathcal{L}_0^{-1}\mathcal{L}$, which is not
symmetric. However, since $\mathcal{L}_0$ is SPD,
\begin{equation}
    \label{eq:app similarity}
    \mathcal{L}_0^{-1}\mathcal{L} = \mathcal{L}_0^{-1/2} \left( \mathcal{L}_0^{-1/2}\, \mathcal{L}\, \mathcal{L}_0^{-1/2} \right) \mathcal{L}_0^{1/2},
\end{equation}
so $\mathcal{L}_0^{-1}\mathcal{L}$ is similar to the symmetric matrix
$S = \mathcal{L}_0^{-1/2}\mathcal{L}\,\mathcal{L}_0^{-1/2}$, and both share the same eigenvalues. By the Courant--Fischer
theorem, the extreme eigenvalues of $S$ are the minimum and maximum of $\mathbf{v}^TS\mathbf{v}/\mathbf{v}^T\mathbf{v}$
over $\mathbf{v} \neq 0$, which, after the change of variables $\mathbf{u} = \mathcal{L}_0^{-1/2}\mathbf{v}$, become the
minimum and maximum of the generalized Rayleigh quotient
\begin{equation}
    \label{eq:app rayleigh}
    R(\mathbf{u}) = \frac{\mathbf{u}^T\mathcal{L}\,\mathbf{u}}{\mathbf{u}^T\mathcal{L}_0\,\mathbf{u}}.
\end{equation}
Therefore, if we manage to show that $\lambda_{\min}(K) \le R(\mathbf{u}) \le \lambda_{\max}(K)$ for every
$\mathbf{u} \neq 0$, which is exactly the statement
$\lambda_{\min}(K)\,\mathcal{L}_0 \preceq \mathcal{L}(\rho) \preceq \lambda_{\max}(K)\,\mathcal{L}_0$, all the eigenvalues
of $\mathcal{L}_0^{-1}\mathcal{L}$ lie in $[\lambda_{\min}(K), \lambda_{\max}(K)]$, and its condition number, i.e. the
ratio of its extreme eigenvalues, is at most $\lambda_{\max}(K)/\lambda_{\min}(K)$. The rest of this appendix is devoted
to proving this two-sided bound on $R(\mathbf{u})$.

\subsection{Making the diffusion tensor explicit in the operator}
The first step is to rewrite Eq.~\eqref{eq:anisotropic operator} in a form in which the entries of $K$ appear
explicitly. We use two properties of the first difference matrices. On the one hand, with Dirichlet boundary conditions
$\mathcal{C}$ is antisymmetric, $\mathcal{C}^T = -\mathcal{C}$, and so is every $\mathcal{C}_\ell$. On the other hand,
for $\ell \neq m$, $\mathcal{C}_\ell$ and $\mathcal{C}_m$ act on different dimensions and therefore commute. Together,
these give $\mathcal{C}_\ell^T\mathcal{C}_m = -\mathcal{C}_\ell\mathcal{C}_m$ for $\ell \neq m$, so that, since
$K_{\ell\ell} = 1$ and $K_{\ell m} = \rho$ for $\ell \neq m$,
\begin{equation}
    \label{eq:app K sum}
    \sum_{\ell,m=1}^N K_{\ell m}\, \mathcal{C}_\ell^T\mathcal{C}_m
    = \sum_{\ell=1}^N \mathcal{C}_\ell^T\mathcal{C}_\ell + \rho \sum_{\ell \neq m} \mathcal{C}_\ell^T\mathcal{C}_m
    = \sum_{\ell=1}^N \mathcal{C}_\ell^T\mathcal{C}_\ell - 2\rho \sum_{\ell < m} \mathcal{C}_\ell\mathcal{C}_m.
\end{equation}
The off-diagonal part of $K$ thus reproduces exactly the mixed-derivative term of $\mathcal{L}(\rho)$, while its diagonal
part produces the extra terms $\mathcal{C}_\ell^T\mathcal{C}_\ell$, which are not present in Eq.~\eqref{eq:anisotropic operator}. Subtracting them from the second difference terms yields the exact identity
\begin{equation}
    \label{eq:app rewrite}
    \mathcal{L}(\rho) = \sum_{\ell=1}^N \left( \mathcal{A}_\ell - \mathcal{C}_\ell^T\mathcal{C}_\ell \right)
    + \sum_{\ell,m=1}^N K_{\ell m}\, \mathcal{C}_\ell^T\mathcal{C}_m.
\end{equation}

\subsection{Splitting the energy into two contributions}
Let us now fix any non-zero grid vector $\mathbf{u}$ and denote its entries by $u_{\mathbf{i}}$, with
$\mathbf{i} = (i_1, \dots, i_N)$ a multi-index over the grid points. We collect at each grid point the discrete gradient
obtained with centered differences,
\begin{equation}
    \label{eq:app gradient}
    \mathbf{g}_{\mathbf{i}} = \left( (\mathcal{C}_1\mathbf{u})_{\mathbf{i}}, \dots, (\mathcal{C}_N\mathbf{u})_{\mathbf{i}} \right)^T \in \mathbb{R}^N,
\end{equation}
and define the two non-negative numbers
\begin{equation}
    \label{eq:app a and c}
    a = \mathbf{u}^T\mathcal{L}_0\,\mathbf{u} = \sum_{\ell=1}^N \mathbf{u}^T\mathcal{A}_\ell\,\mathbf{u}, \qquad
    c = \sum_{\ell=1}^N \|\mathcal{C}_\ell\mathbf{u}\|^2 = \sum_{\mathbf{i}} |\mathbf{g}_{\mathbf{i}}|^2,
\end{equation}
that is, the energy of $\mathbf{u}$ with respect to the Laplacian and the squared norm of its centered discrete gradient,
where $a > 0$ since $\mathcal{L}_0$ is SPD. Since
$\mathbf{u}^T\mathcal{C}_\ell^T\mathcal{C}_m\mathbf{u} = (\mathcal{C}_\ell\mathbf{u})^T(\mathcal{C}_m\mathbf{u})
= \sum_{\mathbf{i}} (\mathcal{C}_\ell\mathbf{u})_{\mathbf{i}}(\mathcal{C}_m\mathbf{u})_{\mathbf{i}}$, the energy of
$\mathbf{u}$ with respect to Eq.~\eqref{eq:app rewrite} reads
\begin{equation}
    \label{eq:app energy}
    \mathbf{u}^T\mathcal{L}(\rho)\,\mathbf{u} = (a - c) + \sum_{\mathbf{i}} \mathbf{g}_{\mathbf{i}}^T K\, \mathbf{g}_{\mathbf{i}}.
\end{equation}
The first contribution corresponds to an isotropic diffusion with unit coefficient, while the second one is the only
place where the anisotropy enters, through the quadratic form of $K$ evaluated on the discrete gradient at every grid
point.

\subsection{Bounding the anisotropic contribution}
For any vector $\mathbf{g} \in \mathbb{R}^N$, the Rayleigh quotient of the SPD matrix $K$ satisfies
$\lambda_{\min}(K)\,|\mathbf{g}|^2 \le \mathbf{g}^TK\mathbf{g} \le \lambda_{\max}(K)\,|\mathbf{g}|^2$. Applying this
inequality at every grid point and summing, the second term of Eq.~\eqref{eq:app energy} is bounded as
\begin{equation}
    \label{eq:app K bound}
    \lambda_{\min}(K)\, c \le \sum_{\mathbf{i}} \mathbf{g}_{\mathbf{i}}^T K\, \mathbf{g}_{\mathbf{i}} \le \lambda_{\max}(K)\, c.
\end{equation}
Hence, it can be written as $\theta\, c$ for some $\theta \in [\lambda_{\min}(K), \lambda_{\max}(K)]$ that depends on
$\mathbf{u}$ (if $c = 0$, any such $\theta$ is valid).

\subsection{Comparing centered and second differences}
The last ingredient is the inequality $0 \le c \le a$, which follows from the one-dimensional property
$0 \preceq \mathcal{C}^T\mathcal{C} \preceq \mathcal{A}$. The lower bound is immediate, since
$\mathbf{u}^T\mathcal{C}^T\mathcal{C}\,\mathbf{u} = \|\mathcal{C}\mathbf{u}\|^2 \ge 0$. For the upper bound, let us
introduce the forward differences $\delta_i= (u_{i+1} - u_i)/h$, for $i = 0, \dots, n$, of a one-dimensional vector
$\mathbf{u}$ with $u_0 = u_{n+1} = 0$. Expanding Eq.~\eqref{eq:app stencils}, one readily checks that
\begin{equation}
    \label{eq:app 1d energies}
    \mathbf{u}^T\mathcal{A}\,\mathbf{u} = \sum_{i=0}^n \delta_i^2, \qquad
    \|\mathcal{C}\mathbf{u}\|^2 = \sum_{i=1}^n \left( \frac{\delta_{i-1} + \delta_{i}}{2} \right)^2,
\end{equation}
that is, the centered difference at a grid point is the average of the two forward differences around it. Since the
square of an average never exceeds the average of the squares, $\left((\delta_{i-1} + \delta_{i})/2\right)^2 \le (\delta_{i-1}^2 +
\delta_{i}^2)/2$, and summing over $i$,
\begin{equation}
    \label{eq:app 1d bound}
    \|\mathcal{C}\mathbf{u}\|^2 \le \frac{1}{2}\sum_{i=1}^n \delta_{i-1}^2 + \frac{1}{2}\sum_{i=1}^n \delta_{i}^2 \le \sum_{i=0}^n \delta_{i}^2 = \mathbf{u}^T\mathcal{A}\,\mathbf{u},
\end{equation}
which proves $\mathcal{C}^T\mathcal{C} \preceq \mathcal{A}$. This inequality carries over to every dimension, as
$\mathcal{A}_\ell - \mathcal{C}_\ell^T\mathcal{C}_\ell$ is the Kronecker product of the positive semidefinite matrix
$\mathcal{A} - \mathcal{C}^T\mathcal{C}$ with identities, and is therefore positive semidefinite itself. Summing over
$\ell$ in Eq.~\eqref{eq:app a and c}, we obtain $0 \le c \le a$, which also shows that the first term of Eq.~\eqref{eq:app energy} is non-negative.

\subsection{The two-sided bound}
We can now combine all the previous steps. Dividing Eq.~\eqref{eq:app energy} by $a$ and introducing $t = c/a$, the
Rayleigh quotient of Eq.~\eqref{eq:app rayleigh} takes the form
\begin{equation}
    \label{eq:app average}
    R(\mathbf{u}) = (1 - t) \cdot 1 + t \cdot \theta, \qquad t \in [0, 1], \quad \theta \in [\lambda_{\min}(K), \lambda_{\max}(K)].
\end{equation}
This is the point at which the inequality $0 \le c \le a$ plays its role: it guarantees that $t \in [0, 1]$, so that
$R(\mathbf{u})$ is a weighted average of the numbers $1$ and $\theta$, and therefore lies between them. It only remains to
note that $1 \in [\lambda_{\min}(K), \lambda_{\max}(K)]$ as well, since all the diagonal entries of $K$ are equal to one,
and the diagonal entries of a symmetric matrix are themselves Rayleigh quotients,
$K_{\ell\ell} = \mathbf{e}_\ell^TK\mathbf{e}_\ell$. Both numbers being averaged lie in the interval
$[\lambda_{\min}(K), \lambda_{\max}(K)]$, and so does their average, which yields
\begin{equation}
    \label{eq:app sandwich}
    \begin{aligned}
        & \lambda_{\min}(K) \le R(\mathbf{u}) \le \lambda_{\max}(K) \quad \forall\, \mathbf{u} \neq 0 \\
        \Longleftrightarrow \quad
        & \lambda_{\min}(K)\,\mathcal{L}_0 \preceq \mathcal{L}(\rho) \preceq \lambda_{\max}(K)\,\mathcal{L}_0.
    \end{aligned}
\end{equation}
In physical terms, part of the energy of $\mathbf{u}$ experiences an isotropic diffusion with unit coefficient and the
rest experiences the anisotropic tensor $K$, and since the unit coefficient lies itself within the spectrum of $K$, the
total energy can never fall below $\lambda_{\min}(K)$ nor exceed $\lambda_{\max}(K)$ times the energy assigned by the
Laplacian.

\subsection{From the two-sided bound to the floor}
As discussed in the preliminaries, Eq.~\eqref{eq:app sandwich} implies that all the eigenvalues of
$\mathcal{L}_0^{-1}\mathcal{L}$ lie in $[\lambda_{\min}(K), \lambda_{\max}(K)]$, so that
\begin{equation}
    \label{eq:app kappa bound}
    \kappa(\mathcal{L}_0^{-1}\mathcal{L}) \le \frac{\lambda_{\max}(K)}{\lambda_{\min}(K)}.
\end{equation}
The eigenvalues of the equicorrelated tensor $K = (1-\rho)\,\mathbb{I}^{(N)} + \rho\,\mathbf{1}\mathbf{1}^T$ are
obtained directly. Along the diagonal direction, $K\mathbf{1} = (1-\rho)\mathbf{1} + \rho N\mathbf{1} =
(1 + (N-1)\rho)\,\mathbf{1}$, while for any $\mathbf{v}$ orthogonal to $\mathbf{1}$, $\mathbf{1}^T\mathbf{v} = 0$ and
$K\mathbf{v} = (1-\rho)\,\mathbf{v}$. For $\rho \ge 0$, the largest eigenvalue is $\lambda_{\max}(K) = 1 + (N-1)\rho$ and
the smallest one is $\lambda_{\min}(K) = 1 - \rho$, and substituting them in Eq.~\eqref{eq:app kappa bound} recovers
Eq.~\eqref{eq:anisotropy floor},
\begin{equation}
    \label{eq:app floor}
    \kappa(\mathcal{L}_0^{-1}\mathcal{L}) \le \frac{1+(N-1)\rho}{1-\rho} =: F(\rho, N).
\end{equation}
For $-1/(N-1) < \rho < 0$ the same argument holds with the roles of the two eigenvalues exchanged, giving the bound
$(1-\rho)/(1+(N-1)\rho)$ instead.

Finally, the derivation also explains why $F$ is only an upper bound, as observed in the discussion of the discrete floor in Sec.~\ref{subsec:anisotropy_experiment}.
Reaching the extremes of Eq.~\eqref{eq:app average} requires $t \to 1$ and $\theta \to \lambda_{\max}(K)$ (or
$\lambda_{\min}(K)$) simultaneously. The former only holds for smooth vectors, for which centered and forward
differences nearly coincide, whereas oscillatory vectors are nearly invisible to $\mathcal{C}$ and have $t \approx 0$.
The latter requires the discrete gradient to be aligned with the corresponding eigenvector of $K$ at every grid point,
which the homogeneous Dirichlet conditions prevent near the boundary. Both conditions can only be approached on fine
grids, where vectors can be simultaneously smooth on the scale of $h$ and strongly varying along a prescribed direction,
which is why the discrete floor remains below $F$ on the coarse benchmark grid.

\end{document}